\documentclass{gtart}

\def\ifplaintex{\expandafter\ifx\csname documentclass\endcsname\relax}

\def\gtp{{\mathsurround=0pt\it $\cal G\mskip-2mu$eometry \&\ 
$\cal T\!\!$opology $\cal P\!$ublications}}  

\def\recd{{\small Received:\qua\receiveddate\ifx\reviseddate\relax
\else\qquad Revised:\qua\reviseddate\fi\par}} 

\def\lognumber#1{\def\thelognumber{#1}}
\def\volumenumber#1{\def\thevolumenumber{#1}}
\def\volumeyear#1{\def\thevolumeyear{#1}}
\def\papernumber#1{\def\thepapernumber{#1}}
\def\pagenumbers#1#2{\def\startpage{#1}\def\finishpage{#2}}
\def\published#1{\def\publishdate{#1}}

\def\received#1{\def\receiveddate{#1}}
\def\revised#1{\def\reviseddate{#1}}
\def\accepted#1{\def\accepteddate{#1}}

\def\asciiauthors#1{\def\theasciiauthors{#1}}

\def\coverauthors#1{\def\thecoverauthors{#1}}
\long\def\asciiabstract#1{\long\def\theasciiabstract{#1}}
\def\asciikeywords#1{\def\theasciikeywords{#1}}

\let\\\par\let\thelognumber\relax\let\thevolumenumber\relax
\let\thepapernumber\relax\let\thevolumeyear\relax\let\startpage\relax
\let\finishpage\relax\let\publishdate\relax\let\receiveddate\relax
\let\reviseddate\relax\let\accepteddate\relax\let\theasciititle\relax
\let\theasciiauthors\relax
\let\theasciiabstract\relax\let\theasciikeywords\relax

\let\thecoverauthors\relax\let\theasciiemail\relax

\ifplaintex
\font\logobig=cmssbx10 scaled 3836
\font\logomed=cmssbx10 scaled 2557
\else
\font\logobig=cmssbx10 scaled 4200
\font\logomed=cmssbx10 scaled 2800
\fi

\long\def\makeagttitle{   
\count0=\startpage
\agt\hfill      
\hbox to 45truept{\vbox to 0pt{\vglue -13truept{\logomed A\kern -.37em{\logobig 
T}\kern -.38em G}\vss}\hss}
\break
{\small Volume \thevolumenumber\ (\thevolumeyear)
\startpage--\finishpage\nl
Published: \publishdate}

\vglue .25truein

{\parskip=0pt\leftskip 0pt plus
1fil\def\\{\par\smallskip}{\Large\bf\thetitle}\par\medskip} \vglue
0.05truein

{\parskip=0pt\leftskip 0pt plus 1fil\def\\{\par}{\sc\theauthors}
\par\medskip}%
 
\vglue 0.03truein 

{\small\leftskip 25truept\rightskip 25truept{\bf Abstract}\stdspace\theabstract

{\bf AMS Classification}\stdspace\theprimaryclass
\ifx\thesecondaryclass\relax\else; \thesecondaryclass\fi\par
{\bf Keywords}\stdspace \thekeywords\par}\vglue 7truept

}   

\ifplaintex
\font\phead=cmsl9 scaled 950
\font\pnum=cmbx10 scaled 913
\font\pfoot=cmsl9 scaled 950
\headline{\vbox to 0pt{\vskip -4.5mm\line{\small\phead\ifnum
\count0=\startpage ISSN 1472-2739 (on-line) 1472-2747 (printed)
\hfill {\pnum\folio}\else\ifodd\count0\def\\{ }%
\ifx\theshorttitle\relax\thetitle\else\theshorttitle\fi\hfill{\pnum\folio}
\else\def\\{ and }{\pnum\folio}\hfill\ifx\theshortauthors\relax\theauthors
\else\theshortauthors\fi\fi\fi}\vss}}
\footline{\vbox to 0pt{\vglue 0mm\line{\small\pfoot\ifnum\count0=\startpage
\copyright\ \gtp\hfill\else
\agt, Volume \thevolumenumber\ (\thevolumeyear)\hfill\fi}\vss}}
\else
\font=cmsl9 scaled 1050
\font\lnum=cmbx10 
\font=cmsl9 scaled 1050
\makeatletter
\def\@oddhead{{\small\lhead\ifnum\count0=\startpage ISSN 1472-2739 
(on-line) 1472-2747 (printed)\hfill {\lnum\number\count0}\else\ifodd\count0
\def\\{ }\ifx\theshorttitle\relax \thetitle \else\theshorttitle\fi\hfill
{\lnum\number\count0}\else\def\\{ and }{\lnum\number\count0}
\hfill\ifx\theshortauthors\relax 
\theauthors\else\theshortauthors\fi\fi\fi}}\def\@evenhead{\@oddhead}
\def\@oddfoot{\small\lfoot\ifnum\count0=\startpage\copyright\ \gtp\hfill\else
\agt, Volume \thevolumenumber\ (\thevolumeyear)\hfill\fi}
\def\@evenfoot{\@oddfoot}
\makeatother
\fi
\let\maketitlepage\makeagttitle
\let\makeshorttitle\maketitlepage
\let\maketitle\maketitlepage

\newwrite\gtoutfile
\long\gdef\makeheadfile{  
{\def\\{, }\def\s{ }
\immediate\openout\gtoutfile head.xxx
\immediate\write\gtoutfile{Proxy-for: \ifx\theasciiauthors\relax
\theauthors\else\theasciiauthors\fi\s<\ifx\theasciiemail\relax\theemail\else\theasciiemail\fi>}
\immediate\write\gtoutfile{\noexpand\\}
\immediate\write\gtoutfile{Authors: \ifx\theasciiauthors\relax
\theauthors\else\theasciiauthors\fi}
{\def\\{ }\immediate\write\gtoutfile{Title: \ifx\theasciititle\relax
\thetitle\else\theasciititle\fi}}
\immediate\write\gtoutfile{Subj-class: GT or SG, GR etc}
\immediate\write\gtoutfile{MSC-class: \theprimaryclass\ifx\thesecondaryclass\relax\else, \thesecondaryclass\fi}
\immediate\write\gtoutfile{Journal-ref: Algebr. Geom. Topol. \thevolumenumber\s
(\thevolumeyear) \startpage-\finishpage}
\immediate\write\gtoutfile{Comments: Published by Algebraic and
Geometric Topology at}
\immediate\write\gtoutfile{\s\s\s  http://www.maths.warwick.ac.uk/agt/AGTVol\thevolumenumber/agt-\thevolumenumber-\thepapernumber.abs.html}
\immediate\write\gtoutfile{\noexpand\\}
\immediate\write\gtoutfile{}
\ifx\theasciiabstract\relax
\immediate\write\gtoutfile{\theabstract}\else
\immediate\write\gtoutfile{\theasciiabstract}\fi
\immediate\write\gtoutfile{}
\immediate\write\gtoutfile{\noexpand\\}
\immediate\write\gtoutfile{}
\immediate\closeout\gtoutfile}}  

\def\maketitlepage{\makeagttitle\makeheadfile}
\let\makeshorttitle\maketitlepage
\let\maketitle\maketitlepage

\lognumber{41}
\volumenumber{3}
\volumeyear{2003}
\papernumber{41}
\published{17 November 2003}
\pagenumbers{1139}{1166}
\received{1 September 2003}
\revised{9 November 2003}
\accepted{9 November2003}

\usepackage{graphicx,rlepsf,amssymb,amsmath,latexsym}
\usepackage[all]{xy}

\theoremstyle{plain}
\newtheorem{theorem}{Theorem}[section]
\newtheorem{lemma}[theorem]{Lemma}

\newtheorem{corollary}[theorem]{Corollary}
\newtheorem{claim}[theorem]{Claim}

\theoremstyle{definition}

\newtheorem*{acknowledgements*}{Acknowledgements}

\newtheorem{example}[theorem]{Example}
\newtheorem{remark}[theorem]{Remark}

\numberwithin{equation}{section}
\numberwithin{figure}{section}

\def \Z{\mathbb Z}
\def \Q{\mathbb Q}
\def \R{\mathbb R}

\def \Hom{{\rm Hom}}
\def \Map{{\rm Map}}
\def \Spin{{\rm Spin}}
\def \Tors{{\rm Tors}}
\def \Quad{{\rm Quad}}
\def \Id{{\rm Id}}
\def \Y{\mathsf{Y}}

\begin{document}

\title[Cohomology rings, Rochlin function and the Goussarov--Habiro theory]
{Cohomology rings, Rochlin function, linking pairing\\
and the Goussarov--Habiro theory of three--manifolds}

\author{Gw\'ena\"el Massuyeau}
\coverauthors{Gw\noexpand\'ena\noexpand\"el Massuyeau}
\asciiauthors{Gwenael Massuyeau}

\address{Institute of Mathematics of the Romanian Academy\\P.O. Box 
1-764, 014700 Bucharest, Romania}

\email{gwenael.massuyeau@imar.ro}

\asciikeywords{3-manifold, surgery equivalence relation, calculus of claspers, spin structure}


\begin{abstract}
We prove that two closed oriented $3$--manifolds have isomorphic quintuplets
(homology, space of spin structures, linking pairing, cohomology rings, Rochlin function)
if, and only if, they belong to the same class of a certain surgery equivalence relation 
introduced by Goussarov and Habiro.
\end{abstract}
\asciiabstract{We prove that two closed oriented 3-manifolds have
isomorphic quintuplets (homology, space of spin structures, linking
pairing, cohomology rings, Rochlin function) if, and only if, they
belong to the same class of a certain surgery equivalence relation
introduced by Goussarov and Habiro.}

\maketitle
\section{Introduction}

Goussarov and Habiro have developed a theory of finite type invariants 
for compact oriented $3$--manifolds \cite{Goussarov, Habiro,GGP}. 
Their theory is based on a new kind of $3$--dimensional 
topological calculus, called \emph{calculus of claspers}.
In strong connection with their finite type invariants, some equivalence relations
have been studied by Goussarov and Habiro. For any integer $k\geq 1$,
the \emph{$Y_k$--equivalence} is the equivalence relation among compact oriented 
$3$--manifolds generated by positive diffeomorphisms and
surgeries along graph claspers of degree $k$. 
The reader will find the precise definition of the $Y_k$--equivalence
in Section \ref{sec:review} and, waiting for this, will be enlightened 
by the following characterization due to Habiro \cite{Habiro}.
Two manifolds $M$ and $M'$ are $Y_k$--equivalent if, and only if, 
there exists a compact oriented connected surface $\Sigma$ in $M$ 
and an element $h$ of the $k$--th lower central series subgroup 
of the Torelli group of $\Sigma$ such that $M'$ is diffeomorphic to the manifold obtained from $M$
by cutting it along $\Sigma$ and re-gluing it using $h$. In particular, we see that
the $Y_k$--equivalence becomes finer and finer as $k$ increases.

Thus, the problem of characterizing the $Y_k$--equivalence relation
in terms of invariants of the manifolds naturally arises.
In the case $k=1$, this problem has been solved for manifolds without boundary. 
Indeed, a result of Matveev \cite{Matveev}, anterior 
to the Goussarov--Habiro theory, can be re-stated as follows: two closed oriented 
$3$--manifolds are $Y_1$--equivalent if and only if they have isomorphic pairs (homology, linking pairing).
That problem has also  been given a solution in the case $k=2$ 
for a certain class of manifolds with boundary \cite{MM}.

In some situations, spin structures and, more recently, complex spin structures have proved to be of use
to low--dimensional topologists. It happens that the Goussarov--Habiro theory can be refined
to the settings where the compact oriented $3$--manifolds are equipped with those additional structures.
So, the problem of characterizing the $Y_k$--equivalence makes sense in those refined contexts as well. 
In the case $k=1$ and for manifolds without boundary, Matveev's theorem has been extended to the realm of spin manifolds
and complex spin manifolds in \cite{Mas} and \cite{DM} respectively.\\

\emph{In this paper, we deal with the $Y_2$--equivalence for manifolds without boundary.}
It is known that surgery along a graph clasper of degree $2$
preserves triple cup products, as well as Rochlin invariant. 
Also, according to Habiro \cite{Habiro}, two homology $3$--spheres are $Y_2$--equivalent 
if and only if they have identical Rochlin invariant. We prove that, in general,
two closed oriented $3$--manifolds are $Y_2$--equivalent if and only if they have isomorphic quintuplets 
(homology, space of spin structures, linking pairing, cohomology rings, Rochlin function).
We also consider the spin case and, with less emphasis, the complex spin case.
In order to give a precise statement of the results, let us fix some notation
for those classical invariants.

Let us consider a closed oriented $3$--manifold $M$. 
A \emph{spin structure} on $M$ is a trivialization of its oriented tangent bundle, up to homotopy 
on $M$ deprived of one point. We denote by Spin$(M)$ the set of  spin structures of $M$ which,
by obstruction theory, is an affine space over the $\Z_2$--vector space $H^1(M;\Z_2)$.
The corresponding action of $H^1(M;\Z_2)$ on Spin$(M)$ is denoted by
$$
\Spin(M) \times H^1(M;\Z_2) \longrightarrow \Spin(M), \quad (\sigma,y) \longmapsto \sigma + y.
$$
We recall that the \emph{Rochlin function} of $M$ is the map
$$R_M: \Spin(M) \longrightarrow \Z_{16}$$
which assigns to any spin structure $\sigma$ on $M$ the signature modulo $16$
of a compact oriented $4$--manifold $W$ such that $\partial W=M$ and $\sigma$ extends to $W$.
The \emph{linking pairing} of $M$, denoted by
$$\lambda_M: \Tors\ H_1(M;\Z) \times \Tors\ H_1(M;\Z) \longrightarrow \Q/\Z,$$
is a nondegenerate symmetric bilinear pairing which measures how rationally null-homologous 
knots are homologically linked in $M$. Let $\Quad(\lambda_M)$ be the space of its quadratic functions,
ie, maps $q: \Tors\ H_1(M;\Z) \to \Q/\Z$ satisfying 
$$q(x_1+x_2)-q(x_1)-q(x_2)=\lambda_M(x_1,x_2)$$
for any $x_1,x_2 \in \Tors\ H_1(M;\Z)$. Lannes, Latour, Morgan and Sullivan  \cite{LL,MS}
have defined a map
$$q_M:\Spin(M) \longrightarrow \Quad(\lambda_M)$$ 
which assigns to any spin structure $\sigma$
a \emph{linking quadratic function} $q_{M,\sigma}$.
For any integer $n\geq 0$,
$$u_M^{(n)}: H^1(M;\Z_n) \times H^1(M;\Z_n) \times H^1(M;\Z_n) \longrightarrow \Z_n$$ 
will denote the skew-symmetric trilinear map 
given by the evaluation of triple cup products with coefficients in $\Z_n$
on the fundamental class of $M$. One can verify, using Poincar\'e duality, that
the cohomology rings of $M$ (with coefficients in $\Z_n$, $n\geq 0$) 
are determined by those triple cup product forms and the group $H_1(M;\Z)$.
Finally, if $M'$ is another closed oriented $3$--manifold and if $\psi: H_1(M;\Z) \to H_1(M';\Z)$ 
is a homomorphism, it will be convenient to denote by $\psi^{(n)}: H^1(M';\Z_n) \to H^1(M;\Z_n)$
the homomorphism corresponding to $\Hom(\psi,\Z_n)$ via Kronecker evaluations. 
\begin{theorem}
\label{th:Y2}
Two closed connected oriented $3$--dimensional manifolds $M$ and $M'$ 
are $Y_2$--equivalent if, and only if, 
there exist an isomorphism $\psi: H_1(M;\Z) \to H_1(M';\Z)$ and 
a bijection $\Psi: \Spin(M') \to \Spin(M)$ 
such that the following conditions hold.
\begin{itemize}
\item[\rm(a)] For any $x_1,x_2 \in \Tors\ H_1(M;\Z)$, we have 
$$\lambda_{M'}\left(\psi(x_1),\psi(x_2)\right)=\lambda_M(x_1,x_2)\in \Q/\Z.$$
\item[\rm(b)] For any integer $n\geq 0$ and for any $y'_1,y'_2,y'_3 \in H^1(M';\Z_n)$, we have 
$$\textstyle{u_{M'}^{(n)}(y'_1,y'_2,y'_3) = 
u_{M}^{(n)}\left(\psi^{(n)}(y'_1),\psi^{(n)}(y'_2),\psi^{(n)}(y'_3)\right) \in \Z_n.}$$
\item[\rm(c)] For any $\sigma' \in \Spin(M')$, we have
$$R_{M'}(\sigma') = R_M(\Psi(\sigma')) \in \Z_{16}.$$
\item[\rm(d)] The bijection $\Psi$ is \emph{compatible} with the isomorphism $\psi$ 
in the sense that it is affine over $\psi^{(2)}$ and the following diagram is commutative:
$$
\xymatrix{
\Spin(M) \ar[r]^-{q_M} & \Quad\left(\lambda_M\right)\\
\Spin(M') \ar[u]^-{\Psi} \ar[r]_-{q_{M'}} & 
\Quad\left(\lambda_{M'}\right). \ar[u]_{\psi^*}
}
$$
\end{itemize}
\end{theorem}

\begin{theorem}
\label{th:Y2_spin}
Two closed connected spin $3$--dimensional manifolds $(M,\sigma)$ and $(M',\sigma')$ are
$Y_2$--equivalent if, and only if, there exists an isomorphism $\psi: H_1(M;\Z) \to H_1(M';\Z)$
such that the following conditions hold.
\begin{itemize}
\item[\rm(a)] For any $x\in \Tors\ H_1(M;\Z)$, we have 
$$q_{M',\sigma'}\left(\psi(x)\right)=q_{M,\sigma}(x)\in \Q/\Z.$$
\item[\rm(b)] For any integer $n\geq 0$ and for any $y'_1,y'_2,y'_3 \in H^1(M';\Z_n)$, we have 
$$\textstyle{u_{M'}^{(n)}(y'_1,y'_2,y'_3) = 
u_{M}^{(n)}\left(\psi^{(n)}(y'_1),\psi^{(n)}(y'_2),\psi^{(n)}(y'_3)\right) \in \Z_n.}$$
\item[\rm(c)] For any $y' \in H^1(M';\Z_2)$, we have
$$\textstyle{R_{M'}(\sigma' + y')=R_M\left(\sigma+\psi^{(2)}(y')\right) \in \Z_{16}.}$$
\end{itemize}
\end{theorem}
\noindent
A similar result holds for manifolds equipped with a complex spin structure (see Section \ref{sec:characterization},
Theorem \ref{th:Y2_spinc}). Let us now discuss the relationship between Theorem \ref{th:Y2} and some 
previously known results.\\

Let $\Sigma_{g,1}$ be the compact connected oriented surface of genus $g$ with one boundary component.
\emph{Homology cylinders} over $\Sigma_{g,1}$ are homology cobordisms with an extra homological 
triviality condition \cite{Habiro, Goussarov}. Homology cylinders form a monoid which contains 
the Torelli group of $\Sigma_{g,1}$ as a submonoid. Moreover, the Johnson homomorphisms 
and the Birman--Craggs homomorphisms extend naturally to this monoid.
An analog of Johnson's result on the Abelianization of the Torelli group 
of $\Sigma_{g,1}$ \cite{Johnson1} has been proved by Meilhan and the author for homology cylinders \cite{MM}:  
two homology cylinders over $\Sigma_{g,1}$ are $Y_2$--equivalent
if and only if they are not distinguished by the first Johnson homomorphism
nor the Birman--Craggs homomorphisms.
On the other hand, there is a canonical construction producing from
any homology cylinder $h$ over $\Sigma_{g,1}$ (for instance, an element $h$ of the
Torelli group) a closed oriented $3$--manifold 
with first homology group isomorphic to $\Z^{2g}$. More precisely, 
one glues to the mapping torus of $h$, which is
a $3$--manifold with boundary $\partial \Sigma_{g,1}\times \mathbf{S}^1$, 
the solid torus $\partial \Sigma_{g,1}\times \mathbf{D}^2$ along the boundary.
Since Johnson \cite{Johnson2}, it is known (at least for elements of the Torelli group) that
the first Johnson homomorphism and the many Birman--Craggs homomorphisms correspond, 
through that construction, to the triple cup products form and the Rochlin function respectively. 
This results in a connection between Theorem \ref{th:Y2} and that characterization of the $Y_2$--equivalence
for homology cylinders. As a matter of fact, some constructions and arguments from \cite{MM} will be re-used here.

Also, it is worth comparing Theorem \ref{th:Y2} to a result of Cochran, Gerges and Orr.
They have studied in \cite{CGO} another equivalence relation
among closed oriented $3$--manifolds, namely the \emph{$2$--surgery equivalence}.
A $2$--surgery, defined as the surgery along a null-homologous knot with framing number $\pm 1$,
is the elementary move of the Cochran--Melvin theory of finite type invariants \cite{CM}.
While the $Y_2$--equivalence coincides with the relation ``have isomorphic quintuplets
(homology, space of spin structures, linking form, cohomology rings, Rochlin function)''
between closed oriented $3$--manifolds, the $2$--surgery equivalence is the relation
``have isomorphic triplets (homology, linking form, cohomology rings)''. 
Indeed, it can be verified that the $Y_2$--equivalence is finer than the
$2$--surgery equivalence, but this will not be used here.

Finally, we mention a result of Turaev, to which Theorem \ref{th:Y2} is complementary.
Consider quintuplets $\textstyle{\left(H,S,\lambda,\left(u^{(n)}\right)_{n\geq 0},R\right)}$
formed by a finitely generated Abel\-ian group $H$, an affine space $S$ 
over the $\Z_2$--vector space $\Hom(H,\Z_2)$, a nondegenerate symmetric 
bilinear pairing $\lambda: \Tors\ H \times \Tors\ H \to \Q/\Z$,
skew-symmetric trilinear forms $u^{(n)}: \Hom(H;\Z_n)^3 \to \Z_n$ and a function 
$R:S\to \Z_{16}$. Turaev has found in \cite{Turaev} necessary and sufficient algebraic conditions
on such a quintuplet to be realized, up to isomorphisms, as the quintuplet
$$\textstyle{\left(H_1(M),\Spin(M),\lambda_M,\left(u^{(n)}_M\right)_{n\geq 0}, R_M\right)}$$
of a closed oriented  $3$--manifold $M$.\\

The paper is organized as follows. In Section \ref{sec:review}, 
we briefly review calculus of claspers and its refinement
to spin manifolds. Next, in Section \ref{sec:variation}, we recall or precise
how the classical invariants involved in Theorem \ref{th:Y2} behave under the surgery along a graph clasper.
In Section \ref{sec:surgery_map}, we fix a closed spin $3$--manifold $(M,\sigma)$
and associate to it a certain set of $Y_2$--equivalence classes. We define a surgery map from a certain space of abstract graphs to this quotient set,
and we prove this map to be bijective. In Section \ref{sec:characterization}, we derive from that bijectivity
Theorem \ref{th:Y2_spin} and, next, Theorem \ref{th:Y2}.
We also give the analogous result for closed oriented $3$--manifolds equipped with a complex spin structure.
Last section is an appendix containing a few algebraic lemmas needed
to obtain the above results.\\

In the sequel, unless otherwise specified, all manifolds are assumed 
to be $3$--dimensional smooth compact and oriented, 
and the diffeomorphisms are supposed to preserve the orientations.

\begin{acknowledgements*}
The author has been supported by an EURROMMAT Fellowship
at the Institute of Mathematics of the Romanian Academy.
\end{acknowledgements*}

\section{Review of calculus of claspers}

\label{sec:review}

We begin by recalling basic concepts from calculus of claspers.
The reader is refered to \cite{Habiro, GGP} for details and complete expositions, 
or to the monograph \cite{O}.\\

\subsection{A flash review of calculus of claspers}

Graph claspers can be defined as follows. We start
with a finite trivalent graph $\mathsf{G}$ decomposed as $\mathsf{G}_1 \cup \mathsf{G}_2$,
where $\mathsf{G}_1$ is a unitrivalent subgraph of $\mathsf{G}$ and 
$\mathsf{G}_2$ is a union of looped edges of $\mathsf{G}$. We give $\mathsf{G}$
a thickening\footnote{A \emph{thickening} of a graph 
$\mathsf{G}$ can be defined as a $\Z_2$--bundle over $\mathsf{G}$ with fiber $[-1,1]$.} 
with the property to be trivial on each looped edge of $\mathsf{G}_2$, and we consider an embedding $G$
of this thickened graph into the interior of a manifold $M$. 
Then, $G$ is said to be a \emph{graph clasper} in the manifold $M$.
The \emph{leaves} of $G$ are the framed knots in $M$ corresponding to the thickening of $\mathsf{G}_2$.
The \emph{degree} of $G$ is the internal degree\footnote{The \emph{internal degree} of a unitrivalent
graph $\mathsf{G}$ is the number of its trivalent vertices if it is connected, or
is the minimum of the internal degrees of its connected components otherwise.} of $\mathsf{G}_1$.
We assume this degree to be at least $1$.
\begin{example} Using the above notations,
if $\mathsf{G}_1$ is $Y$--shaped (respectively $H$--shaped), 
then the graph clasper $G$ is called a \emph{$Y$--graph} (respectively a \emph{$H$--graph}).
Actually, $Y$--graphs play a specific role in the theory. A $Y$--graph and a $H$--graph have been depicted
on Figure \ref{fig:first_ex} before embedding in a manifold $M$. On these diagrams, 
the bold vertices correspond to the trivalent vertices of $\mathsf{G}_1$ and, 
as everywhere else in the
sequel, the graphs are thickened by the ``blackboard framing convention''.
\begin{figure}[ht!]
\centerline{\small \epsfxsize 3.5truein \epsfbox{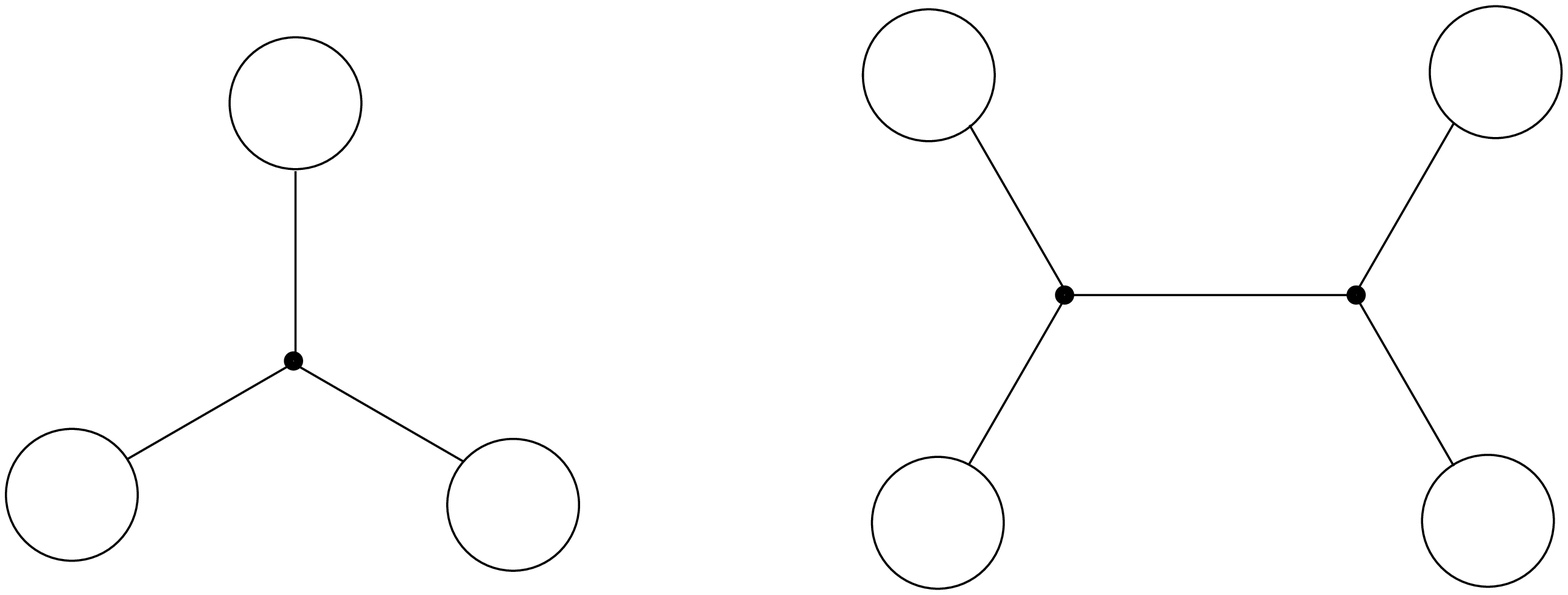}}
\caption{A $Y$--graph and a $H$--graph}
\label{fig:first_ex}
\end{figure} 
\end{example}

A graph clasper carries surgery instructions to modify the manifold where it is embedded.
Surgery along a graph clasper is defined in the following way.

\begin{figure}[ht!]
\centerline{\small \epsfxsize 2.2truein \epsfbox{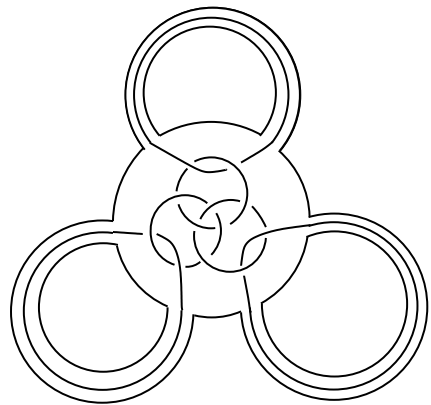}}
\caption{The framed link $B\subset \hbox{N}(G)$}
\label{fig:B}
\end{figure}

First of all, we consider the particular case when $G$ is a $Y$--graph in a manifold $M$.
Let N$(G)$ be the regular neighborhood of $G$ in $M$, which is a genus $3$ handlebody.
The manifold obtained from $M$ by \emph{surgery along} $G$ is denoted by and defined as
$$
M_G:= M\setminus \hbox{int}\left(\hbox{N}(G)\right) \cup_{\partial} \hbox{N}(G)_B
$$
where N$(G)_B$ is N$(G)$ surgered along the six--component framed
link $B$ drawn on Figure \ref{fig:B}.

\begin{figure}[ht!]
\centerline{\small \epsfxsize 4.5truein \epsfbox{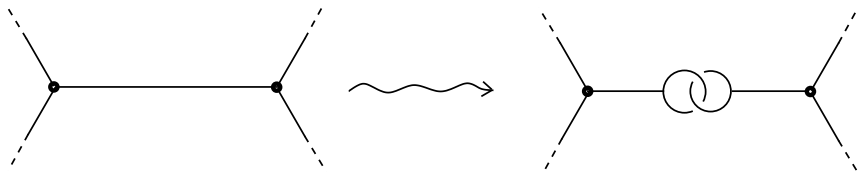}}
\caption{Splitting of a graph clasper to a disjoint union of $Y$--graphs}
\label{fig:splitting}
\end{figure}

Next, we consider the general case when $G$ is a graph clasper in $M$ of arb\-itrary degree $k$. 
By applying the rule illustrated on Figure \ref{fig:splitting}, as many times as necessary, 
$G$ can be transformed to a disjoint union $Y(G)$ of $k$ $Y$--graphs in $M$.
The manifold obtained from $M$ by \emph{surgery along} $G$, denoted by $M_G$,
is the manifold $M$ surgered along each component of $Y(G)$. $M_G$ is also said to be obtained from $M$
by a \emph{$Y_k$--move}. 

The \emph{$Y_k$--equivalence}, mentioned in the introduction, 
is defined to be the equivalence relation among manifolds 
generated by $Y_k$--moves and diffeomorphisms.

\begin{example}
It follows from the definitions that the $Y_1$--equivalence and the $Y_2$--equivalence are generated
by surgeries along $Y$--graphs and $H$--graphs respectively, and diffeomorphisms.
\end{example}

Finally, let us give an idea of what ``calculus of claspers'' is.
Let $G_1$ and $G_2$ be graph claspers in a manifold $M$. 
They are said to be \emph{equivalent}, which we denote by $G_1 \sim G_2$, 
if there exists an embedded handlebody $H$ in $M$ 
whose interior contains both $G_1$ and $G_2$,
and if there exists a diffeomorphism $\tilde{f}:H_{G_1}\to H_{G_2}$
which restricts to the identity on the boundaries
$\partial H_{G_1}\cong \partial H$ and $\partial H_{G_2}\cong \partial H$.
Thanks to the canonical identifications 
$M_{G_i}\cong \left(M\setminus \hbox{int}(H) \right) \cup_\partial H_{G_i}$, 
$\tilde{f}$ induces a diffeomorphism $f: M_{G_1}\to M_{G_2}$.
Hence, $G_1\sim G_2$ implies that $M_{G_1}\cong M_{G_2}$.
The \emph{calculus of claspers} is a corpus of calculi rules which state
equivalence of claspers. Thus, the calculus of claspers allows one to prove diffeomorphisms between manifolds.
\begin{example}
\label{ex:enroulee}
Figure \ref{fig:enroulee} illustrates one of the Goussarov--Habiro moves, which
says that any $H$--graph is equivalent in its regular neighborhood 
(a genus $4$ handlebody) to a $Y$--graph with a null-homologous leaf. 
In particular, the $Y_{k+1}$--equivalence relation is finer than the $Y_k$--equivalence. 
\begin{figure}[h!]
\centerline{\small \epsfxsize 4truein \epsfbox{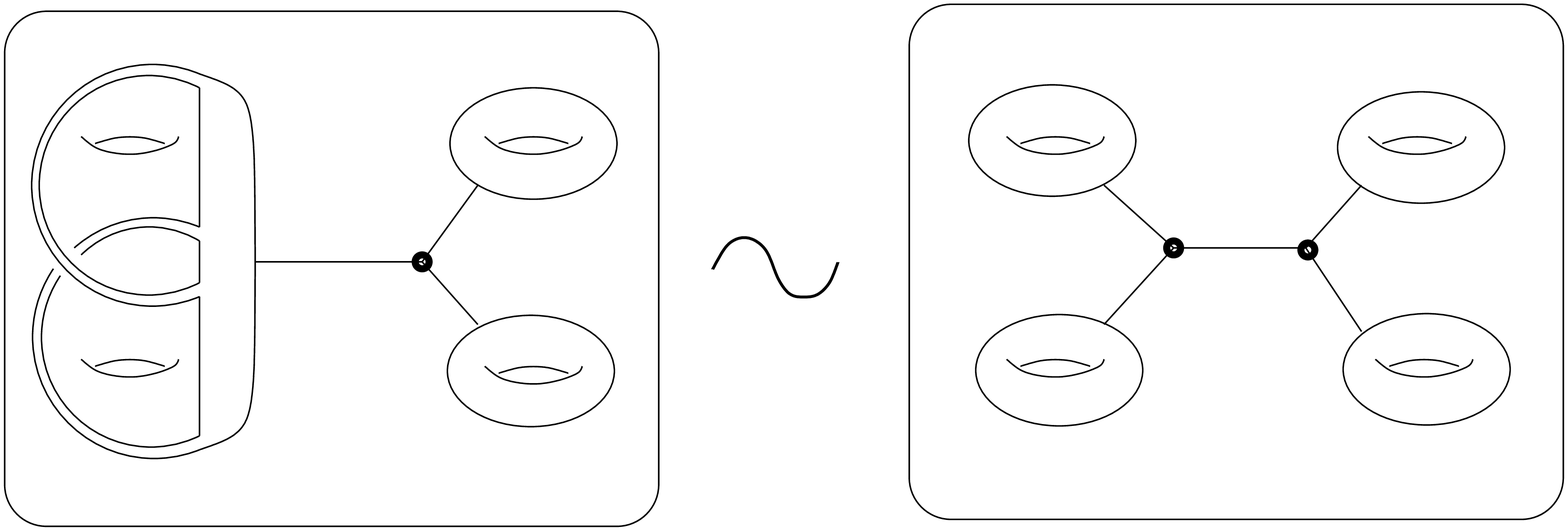}}
\caption{An example of equivalence between graph claspers}
\label{fig:enroulee}
\end{figure}
\end{example}

\subsection{Calculus of claspers for spin manifolds}

\label{subsec:spin}

The most important property of a $Y_k$--move $M \leadsto M_G$ is certainly to preserve 
homology. There is a canonical isomorphism 
$\Phi_G: H_1(M;\Z) \to H_1(M_G;\Z)$, whose existence follows from the fact that
the surgery along a $Y$--graph can be realized by cutting a genus three handlebody and gluing it back
using a diffeomorphism of its boundary which acts trivially in homology \cite{Matveev}.
If $H$ is an embedded handlebody in $M$ whose interior contains $G$, 
$\Phi_G$ is the only map making the diagram
\begin{equation}
\label{diag:Phi}
\xymatrix{
&H_1(M;\Z)\ar@{->}[dd]^{\Phi_G}_\simeq\\
H_1\left(M\setminus H;\Z \right) \ar@{->>}[ru]^{\hbox{\footnotesize incl}_*}
\ar@{->>}[rd]_{\hbox{\footnotesize incl}_*} & \\
&H_1\left(M_G;\Z\right)
}
\end{equation}
commute, where the oblique arrows are induced by inclusions and are surjective.

Furthermore, a $Y_k$--move $M\leadsto M_G$ preserves the space of spin structures.
There exists a canonical bijection $\Theta_G: \Spin(M) \to \Spin(M_G)$, which we shall
denote by $\sigma \mapsto \sigma_G$. This map has been defined in \cite{Mas} for $G$ a $Y$--graph,
the general case can be reduced to this special case by definition of a $Y_k$--move.
If $H$ is a handlebody as above, $\Theta_G$ is the only map making the diagram
\begin{equation}
\label{diag:Theta}
\xymatrix{
\Spin(M) \ar[dd]_{\Theta_G}^\simeq \ \ \ \ar@{>->}[rd]^{\hbox{\footnotesize incl}^*} &\\
& \Spin\left(M\setminus H\right)\\
\Spin(M_G) \ \ \ \ \ar@{>->}[ru]_{\hbox{\footnotesize incl}^*} &
}
\end{equation}
commute, where the oblique arrows are induced by inclusions and are injective.
Let us observe, from diagrams (\ref{diag:Phi}) and (\ref{diag:Theta}), that the bijection 
$\Theta_G$ is affine over the inverse of the isomorphism 
${\Phi_G}^{(2)}: H^1(M_G;\Z_2) \to H^1(M;\Z_2)$ induced by $\Phi_G$.\\

If $G$ is a degree $k$ graph clasper in a manifold $M$ and if $\sigma$ is a spin structure on $M$,
the spin manifold $(M_G,\sigma_G)$ is said to be obtained from
the spin manifold $(M,\sigma)$ by \emph{surgery along} $G$, or, by a \emph{$Y_k$--move}. 
The \emph{$Y_k$--equivalence} among spin manifolds is the equivalence relation generated by
such $Y_k$--moves and spin diffeomorphisms. 
Next lemma says that the calculus of claspers extends to the context of manifolds
equipped with a spin structure.
\begin{lemma}
\label{lem:calculus_spin}
Let $(M,\sigma)$ be a spin manifold. If $G_1$ and $G_2$ are equivalent 
graph claspers in $M$, then the spin manifolds $\left(M_{G_1},\sigma_{G_1}\right)$
and $\left(M_{G_2},\sigma_{G_2}\right)$ are spin diffeomorphic.
\end{lemma}
\begin{proof}
Let $H$ be an embedded handlebody in $M$ whose interior contains $G_1$ and $G_2$,
and let $\tilde{f}:H_{G_1} \to H_{G_2}$ be a diffeomorphism which restricts to the 
identity on the boundaries. Let $f: M_{G_1} \to M_{G_2}$ be the diffeomorphism
induced by $\tilde{f}$. Then, according to (\ref{diag:Theta}),
we have $\Theta_{G_2}=f_* \circ \Theta_{G_1}$. So, $f$ sends 
$\sigma_{G_1}$ to $\sigma_{G_2}$.
\end{proof}
\begin{example}
As in Example \ref{ex:enroulee}, we observe that the $Y_{k+1}$--equivalence is finer
than the $Y_k$--equivalence in the context of spin manifolds too.
\end{example}

\section{Some invariants and surgery along a graph clasper}

\label{sec:variation}

From now on, we restrict ourselves to closed manifolds and, in this section, we describe
how their invariants that are involved in Theorem \ref{th:Y2} behave under the surgery
along a graph clasper.

\subsection{Linking pairing and surgery along a graph clasper}

A theorem of Matveev says that two closed manifolds are $Y_1$--equivalent if and only if they have isomorphic
pairs (homology, linking pairing) \cite{Matveev}. In the spin case,  we have the following refinement of  Matveev's theorem.
\begin{theorem}[See \cite{Mas}]
\label{th:Matveev_spin}
Two closed connected spin manifolds $(M,\sigma)$ and $(M',\sigma')$ are $Y_1$--equivalent
if, and only if, there exists an isomorphism $\psi: H_1(M;\Z) \to H_1(M';\Z)$ such that
\begin{equation}
\label{eq:quadratic}
\forall x\in \Tors\ H_1(M;\Z), \quad q_{M',\sigma'}\left(\psi(x)\right)=q_{M,\sigma}(x)\in \Q/\Z.
\end{equation}
\emph{More precisely}, any isomorphism $\psi: H_1(M;\Z) \to H_1(M';\Z)$ satisfying (\ref{eq:quadratic})
can be realized by a sequence of $Y_1$--moves 
and spin diffeomorphisms from $(M,\sigma)$ to $(M',\sigma')$. 
\end{theorem}
\noindent
Let us comment that characterization.
For any graph clasper $G$ in a closed manifold $M$, we have that
\begin{equation}
\label{eq:Y_linking}
\forall \sigma \in \Spin(M),  \forall x\in \Tors\ H_1(M;\Z),\quad
q_{M_G,\sigma_G} \left(\Phi_G \left(x\right)\right)=q_{M,\sigma}(x).
\end{equation}
This implies the necessary condition in Theorem \ref{th:Matveev_spin}. Reciprocally, 
given an isomorphism $\psi:H_1(M;\Z)\to H_1(M';\Z)$ satisfying (\ref{eq:quadratic}), 
there exists a sequence of $Y_1$--moves and spin diffeomorphisms
$$
(M,\sigma)=(M_0,\sigma_0)\leadsto(M_1,\sigma_1)\leadsto (M_2,\sigma_2) 
\leadsto \cdots\leadsto (M_n,\sigma_n)=(M',\sigma')
$$
such that $\psi=\psi_n\circ \cdots \circ \psi_1$, where
$$
\psi_i=\left\{\begin{array}{ll}
(f_i)_* & \textrm{if } \left(M_{i-1},\sigma_{i-1}\right) \leadsto \left(M_{i},\sigma_{i}\right) 
\textrm{ is a spin diffeomorphism } f_i,\\
\Phi_{G_i}   & \textrm{if } \left(M_{i-1},\sigma_{i-1}\right) \leadsto \left(M_{i},\sigma_{i}\right)
\textrm{ is the surgery along a  $Y$--graph } G_i.
\end{array}\right.
$$
This is what the second statement\footnote{In fact, 
this realization property does not appear explicitely in \cite{Mas} 
but it can be verified from the proof of \cite[Theorem 1]{Mas}.
One of the key ingredients, there, is an algebraic result, due to Durfee and Wall, according to which
two even symmetric bilinear lattices $A$ and $B$ produce isomorphic quadratic functions
$q_A$ and $q_B$ if and only if they are stably equivalent. The point is that, as can be verified from 
\cite[Theorem]{Wall}, any given isomorphism between $q_A$ and $q_B$ can be lifted to a
stable equivalence between $A$ and $B$.} of Theorem \ref{th:Matveev_spin} means.

\subsection{Triple cup products and surgery along a graph clasper}

\label{subsec:Y_triple-cup}

In contrast with the linking quadratic functions, the cohomology rings can be modified
by the surgery along a graph clasper.

\begin{figure}[ht!]
\centerline{\relabelbox \small
\epsfxsize 2truein \epsfbox{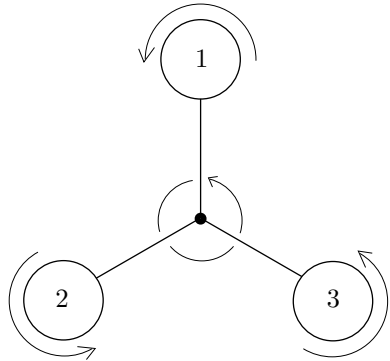}
\adjustrelabel <-0.05cm,0cm> {1}{$1$}
\adjustrelabel <-0.05cm,0cm> {2}{$2$}
\adjustrelabel <-0.05cm,0cm> {3}{$3$}
\endrelabelbox}
\caption{Orientation of each leaf induced by the (cyclic) ordering of the leaves}
\label{fig:orientations}
\end{figure}

\begin{figure}[ht!]
\centerline{\relabelbox \small
\epsfxsize 4truein \epsfbox{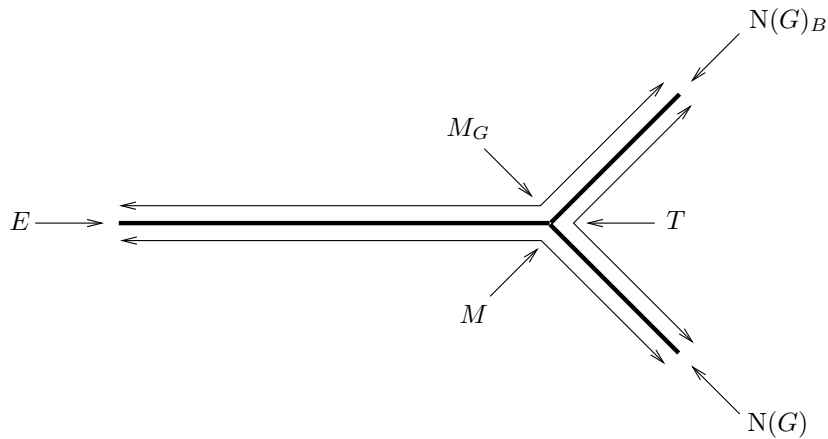}
\adjustrelabel <-0.2cm,-0cm> {E}{$E$}
\adjustrelabel <-0.2cm,-0.15cm> {M}{$M$}
\adjustrelabel <-0.2cm,0.1cm> {MG}{$M_G$}
\relabel {N(G)}{$\hbox{N}(G)$}
\relabel {N(G)B}{$\hbox{N}(G)_B$}
\relabel {torus}{$T$}
\endrelabelbox}
\caption{The singular manifold $N$}
\label{fig:singular}
\end{figure}

Let  $M$ be a closed manifold. For any integer $n\geq 0$, we consider the bilinear pairing
$$\langle-,-\rangle^{(n)}: \Lambda^3 H^1(M;\Z_n)\times \Lambda^3 H_1(M;\Z)  \longrightarrow \Z_n$$ defined by 
\begin{equation}
\label{eq:pairing1}
\langle y_1\wedge y_2 \wedge y_3,x_1\wedge x_2 \wedge x_3\rangle^{(n)}
:= \sum_{\sigma \in {\rm S}_3} \varepsilon(\sigma)\cdot
\prod_{i=1}^3 \langle y_{\sigma(i)},x_i\rangle.
\end{equation}
A $Y$--graph $G$ in $M$ defines an element of $\Lambda^3 H_1(M;\Z)$ 
in the following way. Order the leaves of $G$
and denote them by $L_1,L_2,L_3$ accordingly: $L_1<L_2<L_3$.
This ordering induces an orientation for each leaf, as shown in Figure \ref{fig:orientations}.
Let $\left[L_i\right]$ be the homology class of the $i$--th oriented leaf. Clearly,
$$
\left[L_1\right]\wedge \left[L_2\right] \wedge \left[L_3\right] \in \Lambda^3 H_1(M;\Z)
$$
only depends on $G$. Recall that ${\Phi_G}^{(n)}:H^1(M_G;\Z_n) \to H^1(M;\Z_n)$ stands
for the isomorphism induced by $\Phi_G$.
\begin{lemma}
\label{lem:Y_triple-cup}
Let $G$ be a $Y$--graph in a closed manifold $M$ whose leaves are ordered, 
denoted by $L_1, L_2, L_3$ accordingly 
and oriented as shown in Figure \ref{fig:orientations}. Then, for any integer $n \geq 0$ 
and $y_1',y_2',y_3' \in H^1(M_G;\Z_n)$, we have that
\begin{eqnarray*}
&&\textstyle{u^{(n)}_{M_G}\left(y_1',y_2',y_3'\right)- u^{(n)}_M\left({\Phi_{G}}^{(n)}(y'_1),
{\Phi_{G}}^{(n)}(y'_2), {\Phi_{G}}^{(n)}(y'_3)\right)}\\
&= &\textstyle{\left\langle {\Phi_{G}}^{(n)}(y_1')\wedge {\Phi_{G}}^{(n)}(y_2') \wedge {\Phi_{G}}^{(n)}(y_3'),
[L_1]\wedge [L_2] \wedge [L_3]\right\rangle^{(n)}\in \Z_n.}
\end{eqnarray*}
\end{lemma}
\proof
Let $E:=M\setminus \hbox{int}\left(\hbox{N}(G)\right)$ be the exterior of the $Y$--graph $G$
and consider the singular manifold 
$$
N:=  E \cup_{\partial} \left(\hbox{N}(G)\ \dot\cup\ \hbox{N}(G)_B\right),
$$
which contains both $M$ and $M_G$ (see Figure \ref{fig:singular}). 

Another submanifold of $N$ is $T:=\left(-\hbox{N}(G)\right)\cup_\partial \hbox{N}(G)_B$, 
which is diffeomorphic to the $3$--torus. The group $H_1(T;\Z)$ is free Abelian with basis $(e_1,e_2,e_3)$,
where $e_i$ denotes the homology class of the leaf $L_i$ in N$(G)\subset T$. 
If $e_i^*\in H^1(T;\Z_n)$ is defined by $\langle e_i^*,e_j\rangle=\delta_{ij} \in \Z_n$ for all $i,j=1,2,3$, then 
the cohomology ring of the $3$--torus is such that
\begin{equation}
\label{eq:torus}
\langle e_1^*\cup e_2^* \cup e_3^* , [T]\rangle = 1 \in \Z_n.
\end{equation}

The inclusions induce isomorphisms between $H^1(N;\Z_n)$ and $H^1(M;\Z_n)$, as well
as between  $H^1(N;\Z_n)$ and $H^1(M_G;\Z_n)$.
Let $z_i\in H^1(N;\Z_n)$ be such that $\hbox{incl}^*(z_i)=y'_i \in H^1(M_G;\Z_n)$.
Then, by definition of ${\Phi_G}^{(n)}$, we deduce that $\hbox{incl}^*(z_i)={\Phi_G}^{(n)}(y'_i) \in H^1(M;\Z_n)$.
So, we obtain that
\begin{eqnarray*}
&& \textstyle{\left\langle y'_1\cup y'_2 \cup y'_3, [M_G]\right\rangle-
\left\langle {\Phi_G}^{(n)}(y'_1)\cup {\Phi_G}^{(n)}(y'_2) \cup {\Phi_G}^{(n)}(y'_3),[M]\right\rangle}\\
&=& \left\langle \hbox{incl}^*(z_1)\cup \hbox{incl}^*(z_2) \cup \hbox{incl}^*(z_3),[M_G]\right\rangle\\
&& -\left\langle \hbox{incl}^*(z_1)\cup \hbox{incl}^*(z_2) \cup \hbox{incl}^*(z_3),[M]\right\rangle\\
&=& \left\langle z_1\cup z_2 \cup z_3, \hbox{incl}_*([M_G])-\hbox{incl}_*([M])\right\rangle\\
&=& \left\langle z_1\cup z_2 \cup z_3, \hbox{incl}_*([T])\right\rangle\\
&=& \left\langle \hbox{incl}^*(z_1)\cup \hbox{incl}^*(z_2) \cup \hbox{incl}^*(z_3), [T]\right\rangle.
\end{eqnarray*}
Since $\textstyle{\left\langle \hbox{incl}^* (z_i), e_j \right\rangle=
 \left\langle {\Phi_G}^{(n)}(y'_i),[L_j]\right\rangle}$, we have that 
$$
\hbox{incl}^*(z_i) = \sum_{j=1}^3 \textstyle{\left\langle {\Phi_G}^{(n)}(y'_i),[L_j]\right\rangle e^*_j \in H^1(T;\Z_n)}.
$$
We conclude from (\ref{eq:torus}) that
\begin{align*}
&\textstyle{u^{(n)}_{M_G}\left(y_1',y_2',y_3'\right)- 
u^{(n)}_M\left({\Phi_{G}}^{(n)}(y'_1),{\Phi_{G}}^{(n)}(y'_2), {\Phi_{G}}^{(n)}(y'_3)\right)}\\
&=\textstyle{ \det\left(\left\langle {\Phi_G}^{(n)}(y'_i),[L_j]\right\rangle\right)_{i,j=1,2,3}}\\
&=\textstyle{\left\langle {\Phi_{G}}^{(n)}(y_1')\wedge {\Phi_{G}}^{(n)}(y_2') \wedge {\Phi_{G}}^{(n)}(y_3'),
[L_1]\wedge [L_2] \wedge [L_3] \right\rangle^{(n)}\in \Z_n.}\tag*{\qed}
\end{align*}

\begin{remark}
Lemma \ref{lem:Y_triple-cup} essentially appears in \cite[Section 4.3]{Turaev} where
``Borromean replacements'' are performed on surgery presentations of the manifolds 
in $\mathbf{S}^3$. Indeed, this operation has been used by Turaev to prove his result, mentioned
in the introduction, on realization of skew-symmetric trilinear forms as triple cup products forms of manifolds.
\end{remark}
\begin{corollary} 
\label{cor:triple-cup}
Let $H$ be a graph clasper in a closed manifold $M$ of degree at least $2$.
Then, for any integer $n \geq 0$ 
and $y_1',y_2',y_3' \in H^1(M_H;\Z_n)$, we have that
$$u^{(n)}_{M_H}\left(y_1',y_2',y_3'\right)=
u^{(n)}_M\left({\Phi_{H}}^{(n)}(y_1'),{\Phi_{H}}^{(n)}(y_2'),{\Phi_{H}}^{(n)}(y_3')\right) \in \Z_n.$$
\end{corollary}
\proof
We can suppose that $H$ is connected. By Example \ref{ex:enroulee}, 
$H$ is equivalent to a $Y$--graph $G$ with a null-homologous leaf. If $f: M_{G} \to M_{H}$ is a diffeomorphism
induced by this equivalence of graph claspers, we have  that $ \Phi_{H}=f_*\circ \Phi_{G}$ by diagram (\ref{diag:Phi}).
Applying Lemma \ref{lem:Y_triple-cup} to $G$, we get
\begin{align*}
u^{(n)}_{M_H}\left(y_1',y_2',y_3'\right)& = 
\textstyle{u^{(n)}_{M_{G}}\left(f^*(y_1'),f^*(y_2'),f^*(y_3')\right)}\\
&=\textstyle{u^{(n)}_M\left({\Phi_{G}}^{(n)}f^*(y_1'),{\Phi_{G}}^{(n)}f^*(y_2'),{\Phi_{G}}^{(n)}f^*(y_3')\right)}\\
&=\textstyle{u^{(n)}_M\left({\Phi_{H}}^{(n)}(y_1'),{\Phi_{H}}^{(n)}(y_2'),{\Phi_{H}}^{(n)}(y_3')\right).}\tag*{\qed}
\end{align*}

\subsection{Rochlin invariant and surgery along a graph clasper}

\label{subsec:Y_Rochlin}

As the cohomology rings, the Rochlin invariant can be changed by the surgery along a graph clasper.

Let $M$ be a closed manifold and let ${\hbox{F}}M$ be its bundle of oriented frames, 
which is a $\hbox{GL}_+(3;\R)$--principal bundle:
$$
\xymatrix{
\hbox{GL}_+(3;\R)\ \ar@{>->}[r] & E({\hbox{F}}M) \ar@{->>}[r]^p & M.
}
$$
Let $s \in H_1\left(E({\hbox{F}}M);\Z\right)$ be the image 
of the generator of $H_1\left(\hbox{GL}_+(3;\R);\Z\right)$, which is isomorphic to $\Z_2$.
In this context, the space of spin structures on $M$ can be re-defined as 
\begin{displaymath}
\Spin(M):=\left\{y\in H^1\left(E({\hbox{F}}M);\Z_2\right),
\ \langle y,s \rangle \neq 0\right\}
\end{displaymath}
and the canonical action of $H^1(M;\Z_2)$ on Spin$(M)$ then writes
\begin{equation}
\label{eq:affine}
\forall y\in H^1(M;\Z_2),\ \forall \sigma \in \hbox{Spin}(M), \
\sigma+y:= \sigma + p^*(y).
\end{equation}
(For equivalences between the various definitions of a spin structure, 
the reader is refered to \cite{Milnor}.)

An element $t_K \in H_1\left(E({\hbox{F}}M);\Z\right)$ can be associated
to any oriented framed knot $K\subset M$ in the following way: 
add to $K$ an extra $(+1)$--twist and, next, consider the homology class of its lift in ${\hbox{F}}M$.
Some elementary properties of the map $K \mapsto t_K$ are listed in \cite[Lemma 2.7]{MM}.
\begin{lemma}[See \cite{MM}]
\label{lem:Y_Rochlin}
Let $G$ be a $Y$--graph in a closed manifold $M$ whose leaves 
are ordered, denoted by $L_1, L_2, L_3$ and oriented.  
Then, for any spin structure $\sigma$ on $M$, we have that
\begin{equation}
\label{eq:Rochlin}
R_{M_G}(\sigma_G)-R_M(\sigma)=8\cdot\prod_{k=1}^3 
\langle\sigma,\left[t_{L_k}\right]\rangle\in \Z_{16},
\end{equation}
where $8\cdot :\Z_2 \to \Z_{16}$ denotes the usual monomorphism of groups.
\end{lemma}
\begin{corollary}
\label{cor:Rochlin}
Let $H$ be a graph clasper in a closed manifold $M$ of degree at least two.
Then, for any $\sigma \in \Spin(M)$, we have that $R_{M_H}(\sigma_H)=R_M(\sigma)\in \Z_{16}.$
\end{corollary}
\begin{proof}
Again, we can suppose that $H$ is connected and, by Example \ref{ex:enroulee}, 
$H$ is equivalent to a $Y$--graph $G$ with a null-homologous leaf. By Lemma \ref{lem:calculus_spin},
$\left(M_H,\sigma_H\right)$ is spin diffeomorphic to $\left(M_{G},\sigma_{G}\right)$, hence 
$R_{M_H}\left(\sigma_H\right)=R_{M_{G}}\left(\sigma_{G}\right)$. It follows from \cite[Lemma 2.7]{MM} that
$t_K=0$ for any null-homologous oriented knot $K$ with $0$--framing. 
So, by Lemma \ref{lem:Y_Rochlin}, we have that $R_{M_{G}}(\sigma_{G})=R_M(\sigma).$
\end{proof}

\section{A surgery map}

\label{sec:surgery_map}

In this section, we fix a closed spin manifold $(M,\sigma)$. We associate to $(M,\sigma)$ a bijective surgery map
from a certain space of abstract graphs to a certain set of $Y_2$--equivalence classes.
This is a refinement of the surgery map defined in \cite[Section 2.3]{MM}.

\subsection{Domain and codomain of the surgery map}

\label{subsec:do_codo}

We are going to consider the triplets 
$$\left(M',\sigma',\psi\right),$$
where $(M',\sigma')$ is a spin manifold and $\psi: H_1(M;\Z) \to H_1(M';\Z)$ is an isomorphism such that
$q_{M',\sigma'} \left( \psi \left(x\right)\right)=q_{M,\sigma}(x)$, for any $x\in \Tors\ H_1(M;\Z)$.
The set of such triplets is denoted by
$$
\mathcal{C}(M,\sigma).
$$
A \emph{diffeomorphism} between two triplets 
$\left(M'_1,\sigma'_1,\psi_1\right)$ and $\left(M'_2,\sigma'_2,\psi_2\right)$ in $\mathcal{C}(M,\sigma)$
is a spin diffeomorphism $f:(M'_1,\sigma'_1) \to (M'_2,\sigma'_2)$ such that
$\psi_2=f_*\circ \psi_1$. There is a notion of \emph{$Y_k$--move}, too, for such triplets:
given $\left(M',\sigma',\psi\right) \in \mathcal{C}(M,\sigma)$ and a graph clasper $G$ of degree $k$ in $M'$, 
equation (\ref{eq:quadratic}) allows us to set
$$
\left(M',\sigma',\psi\right)_G := \left(M'_G,\sigma'_G,\Phi_G\circ \psi\right) \in \mathcal{C}(M,\sigma).
$$
The \emph{$Y_k$--equivalence} in $\mathcal{C}(M,\sigma)$ 
is the equivalence relation in $\mathcal{C}(M,\sigma)$ generated by diffeomorphisms and $Y_k$--moves.
The codomain of the surgery map will be the quotient set
$$
\overline{\mathcal{C}}(M,\sigma):=\mathcal{C}(M,\sigma)/Y_2.
$$

Let us now recall a functor defined in \cite[Section 2.1]{MM}.
Let $\mathcal{A}b$ be the category of Abelian groups.  
An \emph{Abelian group with special element} is a pair $(A,s)$ 
where $A$ is an Abelian group and $s\in A$ is of order at most $2$.   
We denote by $\mathcal{A}b_s$ the category of  Abelian groups with special element
whose morphisms are group homomorphisms respecting the special elements. 
We define a functor 
$$
\mathcal{Y}: \mathcal{A}b_s \longrightarrow  \mathcal{A}b 
$$
in the following way. For an object $(A,s)$ of $\mathcal{A}b_{s}$, $\widetilde{\mathcal{Y}}(A,s)$ is defined 
to be the free Abelian group generated by abstract Y--shaped graphs, 
whose edges are given a cyclic order 
and whose univalent vertices are labelled by $A$. The notation
 \[ \Y[a_1,a_2,a_3] \]
will stand for the Y--shaped graph
whose univalent vertices are colored by $a_1$,
$a_2$ and $a_3\in A$ in accordance with the cyclic order,
so that our notation is invariant under cyclic permutation of the $a_i$'s. 
The Abelian group $\mathcal{Y}(A,s)$ is the quotient of 
$\widetilde{\mathcal{Y}}(A,s)$ by the following 
relations\footnote{An antisymmetry relation is also required in \cite{MM}, but
this relation is in fact a consequence of the slide and multilinearity relations.}:

\begin{tabular}{rcl}  
\textbf{Multilinearity} & : &  
$\Y[a_1+a_1',a_2,a_3] = \Y[a_1,a_2,a_3] +  
\Y[a_1',a_2,a_3],$ \\[0.3cm] \textbf{Slide} & : &$\Y[a_1,a_1,a_2] =  \Y[s,a_1,a_2].$  
\end{tabular}

If $f:(A,s)\to (A',s')$ is a morphism in $\mathcal{A}b_{s}$, 
$\mathcal{Y}(f)$ is the group homomorphism $\mathcal{Y}(A,s) \to \mathcal{Y}(A',s')$ 
defined by $\Y[a_1,a_2,a_3] \mapsto \Y[f(a_1),f(a_2),f(a_3)]$.\\

Going back to the spin manifold $(M,\sigma)$, we consider 
the bundle of oriented frames ${\hbox{F}}M$ of $M$.
The domain of the surgery map will be the space of abstract graphs
$\mathcal{Y}(P_M)$ associated to the Abelian group with special element
$$P_M=\left(H_1\left(E({\hbox{F}}M);\Z\right),s\right).$$
Here, as in Section \ref{subsec:Y_Rochlin}, 
$s$ is the image of the generator of $H_1\left(\hbox{GL}_+(3;\R);\Z\right)$.

\subsection{A surgery map from $\mathcal{Y}\left(P_M\right)$
to $\overline{\mathcal{C}}(M,\sigma)$}

Let us consider an arbitrary element $X$ of $\widetilde{\mathcal{Y}}\left(P_M\right)$ written as
$$
X=\sum_{j=1}^n \varepsilon^{(j)}\cdot \Y\left[x_1^{(j)},x_2^{(j)},x_3^{(j)}\right] \quad
\textrm{where}\ \varepsilon^{(j)}=\pm 1\ \textrm{and}\ x_i^{(j)} \in P_M.
$$
For each $j=1,\dots,n$, \emph{pick} a $Y$--graph $G_X^{(j)}$ in $M$ whose leaves are ordered, denoted by
$L_1^{(j)},L_2^{(j)},L_3^{(j)}$ accordingly, oriented as shown in Figure \ref{fig:orientations}
and such that 
$$
\left\{\begin{array}{ll}
t_{L_1^{(j)}}=x_1^{(j)}, t_{L_2^{(j)}}=x_2^{(j)}, t_{L_3^{(j)}}=x_3^{(j)} & \textrm{if}\ \varepsilon^{(j)}=+1,\\
t_{L_1^{(j)}}=x_2^{(j)}, t_{L_2^{(j)}}=x_1^{(j)}, t_{L_3^{(j)}}=x_3^{(j)} & \textrm{if}\ \varepsilon^{(j)}=-1.
\end{array} \right.
$$
(Here, the class $t_K\in H_1(E({\hbox{F}}M);\Z)$ associated to an oriented framed knot $K$ in $M$ 
has been defined in Section \ref{subsec:Y_Rochlin}.) Lastly,
take $G_X$ to be \emph{a} disjoint union of such $Y$--graphs $G_X^{(1)},\dots,G_X^{(n)}$.
\begin{lemma} 
For any $X\in \widetilde{\mathcal{Y}}\left(P_M\right)$, the $Y_2$--equivalence class of
$$
\left(M,\sigma,\Id\right)_{G_X}= \left(M_{G_X},\sigma_{G_X},\Phi_{G_X}\right)\in \mathcal{C}(M,\sigma)
$$
does not depend on the choice of the graph clasper $G_X$ respecting the above requirements.
Moreover, the induced map $\widetilde{\mathcal{Y}}\left(P_M\right) \to \overline{\mathcal{C}}(M,\sigma)$
factors to a quotient map 
$$\mathfrak{S}:\mathcal{Y}\left(P_M\right) \longrightarrow \overline{\mathcal{C}}(M,\sigma).$$
\end{lemma}
\begin{proof}
The demonstration of the lemma, which relies on calculus of claspers,
is very similar to the one given for homology cylinders in \cite[Theorem 2.11]{MM}, so we omit it.
Let us observe that the fact of taking into account, in the definition of $\mathcal{C}(M,\sigma)$,
spin structures together with identifications between the first homology groups 
does not raise extra problems. Indeed, following the proof of Lemma \ref{lem:calculus_spin},
we see that if $G_1$ and $G_2$ are two equivalent graph claspers in $M$, then the triplets
$\left(M_{G_1},\sigma_{G_1},\Phi_{G_1}\right)$ and $\left(M_{G_2},\sigma_{G_2},\Phi_{G_2}\right)$
are diffeomorphic.
\end{proof}

\subsection{Bijectivity of the surgery map $\mathfrak{S}$}

According to the second statement of Theorem \ref{th:Matveev_spin}, 
any element of $\mathcal{C}(M,\sigma)$ is $Y_1$--equivalent to $(M,\sigma,\Id)$.
Consequently, the surgery map $\mathfrak{S}$ is surjective.
In order to prove that $\mathfrak{S}$ is injective too, we are going to insert it into a commutative square
and, for this, we need to define three other maps.
It will be convenient to simplify the notation as follows: 
$S=\Spin(M)$,  $P=P_M$, $H=H_1(M;\Z)$,
$H^{(n)}= \Hom(H,\Z_n)\simeq H^1(M;\Z_n)$ and 
$H_{(n)}=H\otimes \Z_n \simeq H_1(M;\Z_n)$ for any integer $n\geq 0$.\\

Firstly, there is an application
$\textstyle{\overline{\mathcal{C}}(M,\sigma) \to \Map\left(H^{(n)}\times H^{(n)}\times H^{(n)}, \Z_n\right)}$
sending the class of any $(M',\sigma',\psi) \in \mathcal{C}(M,\sigma)$ to the map with value
$$
\textstyle{u^{(n)}_{M'}\left(y'_1,y'_2,y'_3\right)
-u^{(n)}_M\left(\psi^{(n)}(y_1'),\psi^{(n)}(y_2'),\psi^{(n)}(y_3')\right)}
$$
at $\textstyle{\left(\psi^{(n)}(y_1'),\psi^{(n)}(y_2'),\psi^{(n)}(y_3')\right)}$, for any $y'_1,y'_2,y'_3 \in H^1(M';\Z_n)$. 
This map is well-defined because of Corollary \ref{cor:triple-cup}. Similarly,
according to Corollary \ref{cor:Rochlin}, there exists an application
$\overline{\mathcal{C}}(M,\sigma) \to \Map\left(S, \Z_{16}\right)$
sending the class of any $(M',\sigma',\psi) \in \mathcal{C}(M,\sigma)$ to the map with value
$$
\textstyle{R_{M'}(\sigma'+y') -R_M\left(\sigma+\psi^{(2)}(y')\right)}
$$
at $\sigma + \psi^{(2)}(y')$, for any $y'\in H^1(M';\Z_2)$. 
We set
$$
\mathcal{B}(H,S) := \prod_{n\geq 0} 
\Map \textstyle{\left(H^{(n)}\times H^{(n)}\times H^{(n)}, \Z_n\right)} \times \Map\left(S,\Z_{16}\right)
$$ 
and we define
$$
\mathfrak{E}: \overline{\mathcal{C}}(M,\sigma) \longrightarrow \mathcal{B}(H,S)
$$
to be the product of the above maps.\\

Secondly, we come back to the Abelian group with special element $P$.
We denote by $A\left(S,\Z_2\right)$ the space of $\Z_2$--valued affine functions on $S$.
Let $e:H_1\left(E({\hbox{F}}M);\Z\right) \to A\left(S,\Z_2\right)$ be the homomorphism
sending a homology class $x$ to the map $e(x)$ defined by $\alpha \mapsto \langle \alpha,x \rangle$.
(The function $e(x)$ is affine because of (\ref{eq:affine}).)
There exists also a unique homomorphism $\kappa: A(S,\Z_2) \to H_{(2)}$ such that
$f(\sigma+y)=f(\sigma)+\langle y,\kappa(f)\rangle$ for any affine function $f:S\to \Z_2$ 
and cohomology class $y\in H^{(2)}$. Consider the diagram
$$
\xymatrix{
P \ar[r]^-e \ar[d]_-{p_*} & \left(A(S,\Z_2),\overline{1}\right) \ar[d]^-\kappa\\
(H,0) \ar[r]_-{-\otimes \Z_2} & \left(H_{(2)},0\right),
}
$$
in the category of Abelian groups with special element,
where $\overline{1}$ is the
function defined by $\alpha \mapsto 1$ and $p_*$ is the homomorphism in homology induced
by the bundle projection $p:E({\hbox{F}}M) \to M$. By (\ref{eq:affine}), that diagram is commutative:
in fact, according to \cite[Lemma 2.7]{MM}, this is a pull-back square.
In particular, by functoriality, there is a canonical homomorphism
$$
\mathcal{Y}(P) \longrightarrow
\mathcal{Y}(H,0) \times_{\mathcal{Y}(H_{(2)},0)} \mathcal{Y}\left(A(S,\Z_2),\overline{1}\right)
$$
whose codomain is the pull-back of Abelian groups 
obtained from the homomorphisms $\mathcal{Y}(-\otimes \Z_2)$ and $\mathcal{Y}(\kappa)$.
Observe that the groups $\mathcal{Y}(H,0)$ and $\mathcal{Y}(H_{(2)},0)$ are respectively isomorphic
to $\Lambda^3 H$ and $\Lambda^3 H_{(2)}$ via the maps 
defined by $\Y[x_1,x_2,x_3] \mapsto x_1\wedge x_2 \wedge x_3$.
On the other hand, $\mathcal{Y}\left(A(S,\Z_2),\overline{1}\right)$ is isomorphic to
the space of $\Z_2$--valued cubic functions on $S$, denoted by $C(S,\Z_2)$, via the map defined by
$\Y[f_1,f_2,f_3]\mapsto f_1f_2 f_3$. This is proved in the Appendix (Lemma \ref{lem:cubic}). 
Consequently, there is a canonical homomorphism 
$$
\mathfrak{W}: \mathcal{Y}(P) \longrightarrow \Lambda^3 H\times_{\Lambda^3 H_{(2)}} C(S,\Z_2)
$$
whose codomain is the pull-back of Abelian groups obtained from the appropriate 
homomorphisms $\Lambda^3 H \to \Lambda^3 H_{(2)}$ and $C(S,\Z_2) \to \Lambda^3 H_{(2)}$.
The homomorphism $\mathfrak{W}$ is proved to be bijective in the Appendix (Lemma \ref{lem:tri}).\\

Thirdly, there is a homomorphism 
$$
\mathfrak{N}:\Lambda^3 H\times_{\Lambda^3 H_{(2)}} C(S,\Z_2)\longrightarrow \mathcal{B}(H,S)
$$ 
defined by 
$$(X,f)\longmapsto \textstyle{\left(\left(\langle -,X\rangle^{(n)}\right)_{n\geq 0}, 8\cdot f \right)},$$ 
where $\langle -,-\rangle^{(n)}: \Lambda^3 H^{(n)} \times \Lambda^3 H \to \Z_n$ is the pairing defined at (\ref{eq:pairing1}). 
By Lemma \ref{lem:trivectors} from the Appendix, an element $X$ of $\Lambda^3 H$ 
such that $\langle -,X\rangle^{(n)}=0$ for all $n>0$ must vanish. 
Consequently, the homomorphism $\mathfrak{N}$ is injective.\\

The above discussion can be summed up into the square
\begin{equation}
\label{eq:the_square}
\xymatrix{
\Lambda^3H \times_{\Lambda^3 H_{(2)}} C(S,\Z_2)\ \ar@{>->}[r]^-{\mathfrak{N}} & {\mathcal{B}(H,S)}\\
{\mathcal{Y}(P)} \ar@{->>}[r]_-{\mathfrak{S}} \ar[u]_-{\simeq}^-{\mathfrak{W}} & 
\overline{\mathcal{C}}(M,\sigma) \ar[u]_-{\mathfrak{E}}.
}
\end{equation}
The commutativity of that diagram follows from Lemma \ref{lem:Y_triple-cup} and Lemma \ref{lem:Y_Rochlin}. 
We deduce next lemma, which concludes this section on the surgery map $\mathfrak{S}$.
\begin{lemma}
\label{lem:main}
The surgery map $\mathfrak{S}:\mathcal{Y}(P) \to \overline{\mathcal{C}}(M,\sigma)$ 
is bijective, and the map $\mathfrak{E}: \overline{\mathcal{C}}(M,\sigma) \to \mathcal{B}(H,S)$ is injective.
\end{lemma}

\section{Characterization of the $Y_2$--equivalence relation}

\label{sec:characterization}

In this section, we prove the characterization of the $Y_2$--equivalence
relation for closed manifolds, with or without structure, as announced in the introduction.

\subsection{In the setting of spin manifolds: proof of Theorem \ref{th:Y2_spin}}

We start with the necessary condition. 
If $f:(M,\sigma) \to (M',\sigma')$ is a spin diffeomorphism between two closed spin manifolds, 
then conditions (a), (b) and (c) are obviously satisfied for $\psi=f_*:H_1(M;\Z) \to H_1(M';\Z)$. 
Now, we suppose that $G$ is a degree $2$ graph clasper in $M$, we set
$(M',\sigma')=(M_G,\sigma_G)$ and we take $\psi$ to be $\Phi_G: H_1(M;\Z)\to H_1(M_G;\Z)$.
Condition (a) is satisfied, as recalled at (\ref{eq:Y_linking}), and condition (b) too by Corollary \ref{cor:triple-cup}.
Finally, since the bijection $\Omega_G: \Spin(M) \to \Spin\left(M_G\right)$ is affine over the inverse of ${\Phi_G}^{(2)}$,
condition (c) follows from Corollary \ref{cor:Rochlin}.

To prove the sufficient condition, we consider closed spin manifolds $(M,\sigma)$ and $(M',\sigma')$
together with an isomorphism $\psi:H_1(M;\Z)\to H_1(M';\Z)$ satisfying conditions (a), (b) and (c). Then, by (a),
the triplet $(M',\sigma',\psi)$ belongs to $\mathcal{C}(M,\sigma)$ and, by (b) and (c),
$$
\mathfrak{E}(M',\sigma',\psi)=0=\mathfrak{E}(M,\sigma,\Id) \in \mathcal{B}\left(H_1(M;\Z),\Spin(M)\right).
$$ 
Hence, by Lemma \ref{lem:main}, the triplets $(M',\sigma',\psi)$ and $(M,\sigma,\Id)$
are $Y_2$--equivalent in $\mathcal{C}(M,\sigma)$. 
In particular, the spin manifolds $(M,\sigma)$ and $(M',\sigma')$ are $Y_2$--equivalent.
\begin{remark}
\label{rem:precision1}
We have proved a little more than Theorem \ref{th:Y2_spin}: any isomorphism
$\psi:H_1(M;\Z)\to H_1(M';\Z)$ satisfying (a), (b) and (c) can be realized by a sequence of $Y_2$--moves
and spin diffeomorphisms from $(M,\sigma)$ to $(M',\sigma')$.
\end{remark}

\subsection{In the setting of plain manifolds: proof of Theorem \ref{th:Y2}}

Again, the necessary condition is easily verified from previous results. We prove the sufficient condition
and we consider, for this, closed manifolds $M$ and $M'$ together with an isomorphism 
$\psi: H_1(M;\Z) \to H_1(M';\Z)$ and a bijection $\Psi: \Spin(M') \to \Spin(M)$ satisfying conditions
(a) to (d). We choose a spin structure $\sigma'$ on $M'$ and we set $\sigma:=\Psi(\sigma')$.
By condition (d), we have $q_{M,\sigma}(x)=q_{M',\sigma'}(\psi(x))$ for any $x\in \Tors\ H_1(M;\Z)$.
From (d) and (c), we  deduce that
$$
\textstyle{\forall y'\in H^1(M';\Z_2), \ R_M\left(\sigma+\psi^{(2)}(y')\right) = 
R_M\left(\Psi( \sigma'+y')\right)=R_{M'}(\sigma'+y').}
$$
Thus, Theorem \ref{th:Y2_spin} applies: the spin manifolds $(M,\sigma)$ and $(M',\sigma')$ 
are $Y_2$--equivalent and, a fortiori, the manifolds $M$ and $M'$ are $Y_2$--equivalent.
\begin{remark}
\label{rem:precision2}
According to Remark \ref{rem:precision1}, the above proof allows for a more specific statement of Theorem \ref{th:Y2}: 
any pair $(\psi,\Psi)$, formed by an
isomorphism $\psi:H_1(M;\Z)\to H_1(M';\Z)$ and a bijection $\Psi: \Spin(M') \to \Spin(M)$
satisfying conditions (a) to (d), can be realized by a sequence of $Y_2$--moves
and diffeomorphisms from $M$ to $M'$.
\end{remark}

\subsection{In the setting of  complex spin manifolds}

We have seen in Section \ref{subsec:spin} how calculus of claspers makes sense in the context of spin manifolds.
The same happens for manifolds equipped with a complex spin structure. In this paragraph, we give a characterization
of the $Y_2$--equivalence for complex spin manifolds without boundary. Before that, it is worth recalling the
characterization of the $Y_1$--equivalence in this context. 

For a closed manifold $M$, we denote by $B:H_2(M;\Q/\Z) \to \Tors\ H_1(M;\Z)$ the Bockstein
homomorphism associated to the short exact sequence of coefficients $0 \to \Z \to \Q \to \Q/\Z \to 0$.
We also define 
$$L_M: H_2(M;\Q/\Z) \times H_2(M;\Q/\Z) \longrightarrow \Q/\Z$$
to be the symmetric bilinear pairing $\lambda_M \circ (B\times B)$. 
Any complex spin structure $\alpha$ on $M$ produces a quadratic function $\phi_{M,\alpha}$ over $L_M$.
(See \cite{LW,Gille} in case when the Chern class of $\alpha$ is torsion and \cite{DM} in the general case.)
For instance, if $\alpha$ comes from a spin structure $\sigma$, then the quadratic function $\phi_{M,\alpha}$
is essentially equivalent to the linking quadratic function $q_{M,\sigma}$.

According to \cite{DM}, two closed complex spin manifolds
$(M,\alpha)$ and $(M',\alpha')$ are $Y_1$--equivalent if and only if
there exists an isomorphism $\psi:H_1(M;\Z) \to H_1(M';\Z)$ such that 
$\phi_{M,\alpha}\circ \psi^\sharp=\phi_{M',\alpha'}$, where
$\psi^\sharp:H_2(M';\Q/\Z)\to H_2(M;\Q/\Z)$ is the isomorphism dual to $\psi$ 
by the intersection pairings.

\begin{theorem}
\label{th:Y2_spinc}
Two closed connected complex spin $3$--dimensional manifolds 
$(M,\alpha)$ and $(M',\alpha')$ are $Y_2$--equivalent if, and only if, 
there exists an isomorphism $\psi: H_1(M;\Z) \to H_1(M';\Z)$ and 
a bijection $\Psi: \Spin(M') \to \Spin(M)$ such that the following conditions hold.
\begin{itemize}
\item[\rm(a)] For any $z'\in H_2(M';\Q/\Z)$, we have 
$$\textstyle{\phi_{M',\alpha'}\left(z'\right)=\phi_{M,\alpha}\left(\psi^\sharp(z')\right)\in \Q/\Z.}$$
\item[\rm(b)] For any integer $n\geq 0$ and for any $y'_1,y'_2,y'_3 \in H^1(M';\Z_n)$, we have 
$$\textstyle{u_{M'}^{(n)}(y'_1,y'_2,y'_3) = 
u_{M}^{(n)}\left(\psi^{(n)}(y'_1),\psi^{(n)}(y'_2),\psi^{(n)}(y'_3)\right) \in \Z_n.}$$
\item[\rm(c)] For any $\sigma' \in \Spin(M')$, we have
$$\textstyle{R_{M'}(\sigma')= R_M(\Psi(\sigma')) \in \Z_{16}.}$$
\item[\rm(d)] The bijection $\Psi$ is \emph{compatible} with the isomorphism $\psi$ 
in the sense that it is affine over $\psi^{(2)}$ and the following diagram is commutative:
$$
\xymatrix{
\Spin(M) \ar[r]^-{q_M} & \Quad\left(\lambda_M\right)\\
\Spin(M') \ar[u]^-{\Psi} \ar[r]_-{q_{M'}} & 
\Quad\left(\lambda_{M'}\right). \ar[u]_{\psi^*}
}
$$
\end{itemize}
\end{theorem}
\begin{proof}
The necessary condition is proved from previous results (Corollary \ref{cor:triple-cup}, Corollary \ref{cor:Rochlin},
equation (\ref{eq:Y_linking})) and from the following fact: if $G$ is a 
graph clasper in a closed manifold $M$ and if $\alpha$ is a complex spin structure on $M$, then we have that
\begin{equation}
\label{eq:quadratic_spinc}
\textstyle{\forall z'\in H_2\left(M_G;\Q/\Z\right), \
\phi_{M_G,\alpha_G}(z')=\phi_{M,\alpha}\left({\Phi_G}^\sharp (z') \right) \in \Q/\Z.}
\end{equation}
To show the sufficient condition, we consider closed manifolds 
equipped with a complex spin structure $(M,\alpha)$ and $(M',\alpha')$,
together with bijections $\psi$ and $\Psi$ satisfying conditions (a) to (d).
We denote by $\phi_M: \Spin^c(M) \to \Quad\left(L_M\right)$ the map defined by 
$\alpha \mapsto \phi_{M,\alpha}$: it turns out to be injective \cite{DM}.
By condition (a), we have $L_{M'}=L_M\circ\left(\psi^\sharp \times \psi^\sharp \right)$
or, equivalently, $\lambda_{M}=\lambda_{M'}\circ \left(\psi|_{\Tors}\times \psi|_{\Tors}\right)$.
Therefore, by Theorem \ref{th:Y2} and Remark \ref{rem:precision2}, there exists a sequence
of $Y_2$--moves and diffeomorphisms from $M$ to $M'$ which realizes the isomorphism $\psi$
in homology. This sequence of moves induces an identification $\Psi^c$ between $\Spin^c(M')$
and $\Spin^c(M)$ which, by identity (\ref{eq:quadratic_spinc}), makes the diagram
$$
\xymatrix{
\Spin^c(M) \ar[r]^-{\phi_M} & \Quad\left(L_M\right) \ar[d]^{\left(\psi^\sharp\right)^*}\\
\Spin^c(M') \ar[u]^-{\Psi^c} \ar[r]_-{\phi_{M'}} & 
\Quad\left(L_{M'}\right)
}
$$
commute. In particular, we have $\phi_{M,\Psi^c(\alpha')}\circ \psi^\sharp= 
\phi_{M',\alpha'}=\phi_{M,\alpha}\circ \psi^\sharp$, hence $\alpha=\Psi^c(\alpha')$.
We conclude that the complex spin manifolds $(M,\alpha)$ and $(M',\alpha')$ are $Y_2$--equivalent.
\end{proof}

\section{Appendix}

This section contains the proofs of the algebraic lemmas that have been used in Section \ref{sec:surgery_map}.
Here, we shall use the following convention for any finitely generated Abelian group $A$ and any integer $n>0$.
We denote $A_{(n)}=A\otimes \Z_n$, $A^{(n)}=\Hom\left(A,\Z_n\right)$ and 
$\textstyle{\langle-,-\rangle^{(n)}: \Lambda^3 A^{(n)}  \times \Lambda^3 A \to \Z_n}$
the pairing defined by
\begin{equation}
\label{eq:pairing2}
\langle y_1\wedge y_2 \wedge y_3, x_1\wedge x_2 \wedge x_3\rangle^{(n)}
:= \sum_{\sigma \in {\rm S}_3} \varepsilon(\sigma)\cdot
\prod_{i=1}^3 \langle y_{\sigma(i)},x_i\rangle.
\end{equation}
A \emph{basis} of $A$ is a family of pairs $\left\{(e_i,n_i):i\in I\right\}$
indexed by a finite set $I$, such that $e_i$ is 
an element of $A$ of order\footnote{For an element $e$ 
of an Abelian group $A$, the \emph{order} of $e$
 is the unique integer $n\geq 0$ such that the subgroup generated by $e$
is isomorphic to $\Z_n$.} $n_i\geq 0$ and $A$ is the direct sum of the cyclic subgroups generated by the $e_i$'s. 
The \emph{dual} basis of $A^{(n)}$ is the basis $\left\{(e_i^*,\gcd(n,n_i)):i\in I\right\}$ of $A^{(n)}$, 
defined by $\langle e^{*}_i ,e_j\rangle =\delta_{i,j} n/\gcd(n,n_i) \in \Z_n$.

\subsection{Embedding of trivectors}

This paragraph is aimed at proving the following lemma.
\begin{lemma}
\label{lem:trivectors}
If $H$ is a finitely generated Abelian group, the homomorphism
$$
\Lambda^3 H \longrightarrow \prod_{n> 0} \textstyle{\Hom\left(\Lambda^3 H^{(n)},\Z_n\right)}, \
X \longmapsto \textstyle{\left(\left\langle -,X\right\rangle^{(n)}\right)_n}
$$
is injective. 
\end{lemma}
\begin{proof}
Let $X\in \Lambda^3 H$ be such that
$$
\textstyle{(\mathcal{H}_m)\quad \quad \quad \quad \left\langle -,X\right\rangle^{(m)}
= 0 \in \Hom\left(\Lambda^3 H^{(m)},\Z_m\right)}
$$
for all integers $m > 0$. To show that $X$ must vanish, it suffices to prove that
$$
(\mathcal{A}_m)\quad \quad \quad X\otimes 1=0 \in \left(\Lambda^3 H\right)\otimes \Z_m\simeq \frac{\Lambda^3 H}{m\cdot\Lambda^3 H}
$$
for any integer $m> 0$. Assertion $(\mathcal{A}_m)$ trivially holds for $m=1$ so that it suffices to prove the
following inductive statement.
\begin{claim}
\label{claim:induction}
Let $n > 0$ be an integer. 
If assertion $(\mathcal{A}_n)$ holds, then assertion $(\mathcal{A}_{np})$
holds too for any prime number $p$.
\end{claim}
To prove Claim \ref{claim:induction}, we need a few preliminaries. Choose a basis 
$$
\left\{(e_i,n_i):1\leq i \leq r\right\}
$$ 
of $H$, and let $m$ be an arbitrary positive integer. Then,
$$\{(e_i\wedge e_j\wedge e_k,\gcd(n_i,n_j,n_k)) : 1\leq i < j < k\leq r\}$$ 
is a distinguished basis of $\Lambda^3 H$, while
$$\{((e_i\wedge e_j\wedge e_k)\otimes 1,\gcd(m,n_i,n_j,n_k)) : 1\leq i < j < k\leq r\}$$
is a preferred basis of $\left(\Lambda^{3} H\right)\otimes \Z_m$. 
Furthermore, a distinguished basis of $\Lambda^3 H^{(m)}$ is 
$$
\left\{\left(e_i^* \wedge e_j^* \wedge e_k^*
,\gcd(m,n_i,n_j,n_k)\right):1\leq i < j < k\leq r\right\}
$$
and $\Hom\left(\Lambda^3 H^{(m)},\Z_m\right)$ has the basis
$$
\left\{\left(\left(e_i^* \wedge e_j^* \wedge e_k^*\right)^*,
\gcd(m,n_i,n_j,n_k)\right):1\leq i < j < k\leq r\right\}.
$$
The homomorphism $\Lambda^3 H \to \Hom\left(\Lambda^3 H^{(m)},\Z_m\right)$ 
defined by $Y\mapsto \langle -, Y\rangle^{(m)}$ sends the basis element $e_i\wedge e_j\wedge e_k$ to 
\begin{equation}
\label{eq:non-isomorphism}
\frac{m^2 \gcd(m,n_i,n_j,n_k)}{\gcd(m,n_i)\gcd(m,n_j)\gcd(m,n_k)}
\left(e_i^* \wedge e_j^* \wedge e_k^*\right)^*.
\end{equation}

We suppose that $(\mathcal{A}_n)$ holds, we consider a prime number $p$
and we want to show that $(\mathcal{A}_{np})$ holds. Writing $X$ in the preferred basis of $\Lambda^3 H$, say
$$
X=\sum_{1\leq i<j<k\leq r} x_{ijk}\cdot e_i\wedge e_j\wedge e_k \quad  \ (x_{ijk} \in \Z),
$$
this amounts to prove that 
$$x_{ijk}\equiv 0 \mod \gcd(np,n_i,n_j,n_k).$$
But, from $(\mathcal{A}_n)$, we know that
$$
x_{ijk}\equiv 0 \ \mod \gcd(n,n_i,n_j,n_k)
$$
and, from $(\mathcal{H}_{np})$ together with  (\ref{eq:non-isomorphism}) applied to $m=np$, we know that
$$
\frac{x_{ijk}n^2p^2\gcd(np,n_i,n_j,n_k)}{\gcd(np,n_i)\gcd(np,n_j)\gcd(np,n_k)} 
\equiv 0 \mod \gcd(np,n_i,n_j,n_k).
$$
Therefore, it is enough to prove that the conditions
$$
\left\{ \begin{array}{rcl}
z &\equiv & 0 \mod \gcd(n,n_i,n_j,n_k)\\
z n^2 p^2 &\equiv& 0 \mod \gcd(np,n_i)\gcd(np,n_j)\gcd(np,n_k)
\end{array}\right. 
$$
imply that $ z \equiv 0 \mod \gcd(np,n_i,n_j,n_k)$ for any integer $z$.
But, this can be verified working with the $p$--valuations of $n$, $n_i$, $n_j$, $n_k$ and $z$.
\end{proof}

\subsection{Cubic functions and trivectors}

Let $H$ be a finitely generated Abelian group and let $S$ be a $\Z_2$--affine space over $H^{(2)}$.
We denote by $A(S,\Z_2)$ the space of affine functions $S\to \Z_2$ 
and by $\overline{1} \in A(S,\Z_2)$ the constant function $\sigma \mapsto 1$.
Then, $\left(A(S,\Z_2),\overline{1}\right)$ is an Abelian group with special element 
(in the sense of Section \ref{subsec:do_codo}). The space of cubic functions $S\to \Z_2$, 
ie, functions which are finite sums of triple products of affine functions,  
is denoted by $C(S,\Z_2)$.
\begin{lemma}
\label{lem:cubic}
The homomorphism $\gamma: \mathcal{Y}\left(A(S,\Z_2),\overline{1}\right) \to C(S,\Z_2)$
defined by $\gamma\left(\Y[f_1,f_2,f_3]\right)=f_1f_2f_3$ is an isomorphism.
\end{lemma}
\begin{proof}
The demonstration is similar to that of \cite[Lemma 4.21]{MM}. It is enough to construct an epimorphism
$\epsilon : C(S,\Z_2)\to \mathcal{Y}\left(A(S,\Z_2),\overline{1}\right)$ such that 
$\gamma \circ \epsilon$ is the identity.

We fix a base point $\sigma_0 \in S$ together with a basis $\{(e_i,n_i):1\leq i \leq r\}$ of 
the Abelian group $H$. Let $\overline{e_i}: S\to \Z_2$ be the affine function defined by
$\overline{e_i}(\sigma_0+y):= \langle y, e_i \rangle$ for any $y\in H^{(2)}$. Then, 
$$\left\{\left(\overline{1},2\right)\right\} \cup \left\{\left(\overline{e_i},\gcd(2,n_i)\right):1\leq i\leq r\right\}$$ 
is a basis of $A(S,\Z_2)$, and a basis of $C(S,\Z_2)$ is given by
\begin{eqnarray*}
&&\hphantom{\cup}\left\{\left(\overline{1},2\right)\right\} \cup \left\{\left(\overline{e_i},\gcd(2,n_i)\right),1\leq i \leq r\right\}\\
&&\cup \left\{\left(\overline{e_i}\overline{e_j},\gcd(2,n_i,n_j)\right),1\leq i<j \leq r\right\}\\
&&\cup \left\{\left(\overline{e_i}\overline{e_j}\overline{e_k},\gcd(2,n_i,n_j,n_k)\right),1\leq i<j<k \leq r\right\}.
\end{eqnarray*}
We define a homomorphism $\epsilon : C(S,\Z_2)\to \mathcal{Y}\left(A(S,\Z_2),\overline{1}\right)$ by setting
\begin{eqnarray*}
\epsilon\left(\overline{1}\right) &:= &\Y\left[\overline{1},\overline{1},\overline{1}\right],\\
\epsilon\left(\overline{e_i}\right)&:=& \Y\left[\overline{e_i},\overline{1},\overline{1}\right],\\
\epsilon\left(\overline{e_i}\overline{e_j}\right) &:=& \Y\left[\overline{e_i},\overline{e_j},\overline{1}\right],\\
\epsilon\left(\overline{e_i}\overline{e_j}\overline{e_k}\right)&:=& \Y\left[\overline{e_i},\overline{e_j},\overline{e_k}\right].
\end{eqnarray*}
The homomorphism $\epsilon$ is surjective by the slide and multilinearity relations and, clearly, 
the identity $\gamma\circ \epsilon =\Id$ is satisfied on the basis elements.
\end{proof}

Given a cubic function $f:S\to \Z_2$, one can compute its \emph{formal third derivative} 
$$\hbox{d}^3 f: H^{(2)}\times H^{(2)}\times H^{(2)} \longrightarrow \Z_2$$ 
defined for any $\sigma \in S$ by
$$
\hbox{d}^3 f(y_1,y_2,y_3):= \sum_{(\epsilon_1,\epsilon_2,\epsilon_3) \in \{0,1\}^3} 
f\left(\sigma+\epsilon_1\cdot y_1+\epsilon_2 \cdot y_2 +\epsilon_3 \cdot y_3\right).
$$
It can be verified that $\hbox{d}^3 f$ is multilinear, does not depend
on $\sigma$ (because $f$ is cubic) and is alternate (because $2\cdot H^{(2)}=0$), 
hence a homomorphism $\hbox{d}^3: C(S,\Z_2) \to \Hom\left(\Lambda^3 H^{(2)},\Z_2\right)$.
By the duality pairing between $\Lambda^3 H_{(2)}$ and $\Lambda^{3} H^{(2)}$ defined by equation (\ref{eq:pairing2}),
this homomorphism can be regarded as taking its values in $\Lambda^3 H_{(2)}$.
In the sequel, we consider the pull-back of Abelian groups
\begin{equation}
\label{eq:Johnson_pull-back}
\Lambda^3 H\times_{\Lambda^3 H_{(2)}} C(S,\Z_2)
\end{equation}
defined by $\Lambda^3\left(-\otimes \Z_2\right): \Lambda^3 H\to \Lambda^3H_{(2)}$
and  $\hbox{d}^3:C(S,\Z_2) \to \Lambda^{3} H_{(2)}$.

Let now $P$ be the pull-back of Abelian groups with special element 
$$
P:= \left(H,0\right)\times_{\left(H_{(2)},0\right)} \left(A(S,\Z_2),\overline{1}\right)
$$
induced by the homomorphisms $-\otimes \Z_2: H\to H_{(2)}$ and $\kappa: A(S,\Z_2) \to H_{(2)}$, where $\kappa$ is
defined by $f(\sigma+y)=f(\sigma)+\langle y,\kappa(f)\rangle$ for any $f\in A(S,\Z_2)$, $\sigma \in S$ and $y\in H^{(2)}$. 
By functoriality, there is a canonical homomorphism 
$$\mathcal{Y}(P) \longrightarrow \mathcal{Y}(H,0) \times_{\mathcal{Y}(H_{(2)},0)} \mathcal{Y}(A(S,\Z_2),\overline{1})$$
with values in the pull-back of Abelian groups
obtained from the homomorphisms $\mathcal{Y}(-\otimes \Z_2)$ and $\mathcal{Y}(\kappa)$.

The isomorphisms $\mathcal{Y}(H,0)\simeq \Lambda^3 H$
and $\mathcal{Y}\left(H_{(2)},0\right)\simeq \Lambda^3 H_{(2)}$ 
(defined by $\Y[x_1,x_2,x_3]\mapsto x_1\wedge x_2\wedge x_3$)
together with the isomorphism $\gamma$ of Lemma \ref{lem:cubic} induce an isomorphism
\begin{equation}
\label{eq:iso_pull-backs}  
\mathcal{Y}(H,0) \times_{\mathcal{Y}(H_{(2)},0)} \mathcal{Y}(A(S,\Z_2),\overline{1})\simeq
\Lambda^3 H\times_{\Lambda^3 H_{(2)}} C(S,\Z_2)
\end{equation} 
between the above two pull-backs of Abelian groups. 
Indeed, it can be verified that
$\hbox{d}^3(f_1f_2f_3)=\kappa(f_1)\wedge \kappa(f_2) \wedge \kappa(f_3) \in \Lambda^3 H_{(2)}$
for any $f_1,f_2,f_3 \in A(S,\Z_2)$.
\begin{lemma}
\label{lem:tri}
Let $\mathfrak{W}: \mathcal{Y}(P) \to \Lambda^3 H\times_{\Lambda^3 H_{(2)}} C(S,\Z_2)$
be the canonical homomorphism 
$\mathcal{Y}(P)\to \mathcal{Y}(H,0) \times_{\mathcal{Y}(H_{(2)},0)} \mathcal{Y}(A(S,\Z_2),\overline{1})$
composed with the isomorphism (\ref{eq:iso_pull-backs}). Then, $\mathfrak{W}$ is an isomorphism. 
\end{lemma}
\begin{proof}
As in the proof of Lemma \ref{lem:cubic}, it suffices to construct an epimorphism 
$\epsilon: \Lambda^3 H\times_{\Lambda^3 H_{(2)}} C(S,\Z_2) \to \mathcal{Y}(P)$ 
such that $\mathfrak{W} \circ \epsilon$ is the identity. 
Again, we fix a basis $\left\{(e_i,n_i):1\leq i \leq r\right\}$ of $H$
together with a base point $\sigma_0\in S$, and $\overline{e_i}: S\to \Z_2$ designates the affine function
defined by $\overline{e_i}(\sigma_0+ y):= \langle y, e_i \rangle$ for any $y\in H^{(2)}$.
From the basis of $A(S,\Z_2)$ given in the proof of Lemma \ref{lem:cubic}, we obtain that
$$
\left\{ ( (0,\overline{1}) , 2 ) \right\} \cup \left\{ ( (e_i,\overline{e_i}) , n_i ): 1\leq i \leq r \right\}
$$
is a basis of $P$. From the basis of $C(S,\Z_2)$ given in the proof of Lemma \ref{lem:cubic}
and the basis $\left\{(e_i\wedge e_j\wedge e_k,\gcd(n_i,n_j,n_k)):1\leq  i < j < k \leq r\right\}$
of $\Lambda^3 H$, we construct the following basis of 
$\textstyle{\Lambda^3 H\times_{\Lambda^3 H_{(2)}} C(S,\Z_2)}$:
\begin{eqnarray*}
&&\hphantom{\cup} \left\{ \left( \left(0,\overline{1}\right) , 2 \right) \right\} 
\cup \left\{ ((0,\overline{e_i}),\gcd(2,n_i)): 1\leq i \leq r \right\}\\
&&\cup \left\{ ((0,\overline{e_i} \overline{e_j}),\gcd(2,n_i,n_j)): 1\leq i<j \leq r \right\}\\
&& \cup \left\{ ((e_i\wedge e_j\wedge e_k,\overline{e_i} \overline{e_j} \overline{e_k}),\gcd(n_i,n_j,n_k)): 1\leq i<j<k \leq r \right\}.
\end{eqnarray*}
A homomorphism $\epsilon: \Lambda^3 H\times_{\Lambda^3 H_{(2)}} C(S,\Z_2) \to \mathcal{Y}(P)$ is defined
by giving its values on the basis elements in the following way:
\begin{eqnarray*}
\epsilon\left(0,\overline{1}\right) & := &
\Y\left[\left(0,\overline{1}\right),\left(0,\overline{1}\right),\left(0,\overline{1}\right)\right],\\
\epsilon\left(0,\overline{e_i}\right) & := &
\Y\left[(e_i,\overline{e_i}),\left(0,\overline{1}\right),\left(0,\overline{1}\right)\right],\\
\epsilon\left(0,\overline{e_i}\overline{e_j}\right) & := &
\Y\left[(e_i,\overline{e_i}),(e_j,\overline{e_j}),\left(0,\overline{1}\right)\right],\\
\epsilon\left(e_i\wedge e_j\wedge e_k,\overline{e_i}\overline{e_j}\overline{e_k}\right) 
& := &\Y[(e_i,\overline{e_i}),(e_j,\overline{e_j}),(e_k,\overline{e_k})].\\
\end{eqnarray*}
By the slide and multilinearity relations, this homomorphism $\epsilon$ is surjective, 
and it can be readily verified that $\mathfrak{W} \circ \epsilon (z) =z$ for any of the above basis elements $z$. 
\end{proof}

\Addresses\recd

\end{document}
